# A Quick Estimation of Fréchet Quantizers for a Dynamic Solution to Flood Risk Management Problems

Anna Timonina-Farkas
École Polytechnique Fédérale de Lausanne (EPFL), College of Management of Technology,
Chair of Technology and Operations Management,
EPFL-CDM-TOM, Station 5, CH-1015 Lausanne, Switzerland
International Institute for Management Development (IMD),
Chemin de Bellerive 23, 1003 Lausanne, Switzerland, anna.farkas@epfl.ch, anna.farkas@imd.org

Multi-stage stochastic optimization is a well-known quantitative tool for decision-making under uncertainty. It is broadly used in financial and investment planning, inventory control, and also natural disaster risk management. Theoretical solutions of multi-stage stochastic programs can be found explicitly only in very exceptional cases due to their variational form and interdependency of uncertainty in time. Nevertheless, numerical solutions are often inaccurate, as they rely on Monte-Carlo sampling, which requires the Law of Large Numbers to hold for the approximation quality. In this article, we introduce a new approximation scheme, which computes and groups together stage-wise optimal quantizers of conditional Fréchet distributions for optimal weighting of value functions in the dynamic programming. We consider optimality of scenario quantization methods in the sense of minimal Kantorovich-Wasserstein distance at each stage of the scenario tree. By this, we bound the approximation error with convergence guarantees. We also provide global solution guarantees under convexity and monotonicity conditions on the value function. We apply the developed methods to the governmental budget allocation problem for risk management of flood events in Austria. For this, we propose an extremely efficient way to approximate optimal quantizers for conditional Fréchet distributions. Our approach allows to enhance the overall efficiency of dynamic programming via the use of different parameter estimation methods for different groups of quantizers. The groups are distinguished by a particular risk threshold and are able to differentiate between higher- and lower-impact flood events.

## 1. Introduction and Literature Review

Accurate and efficiently computed decisions under time-varying uncertainties are required in a wide variety of real-world applications (Blanchet [12], Calma et al. [14], Eijgenraam et al. [18], Pflug [40], Pflug and Römisch [41], Powell [47]). Multi-stage stochastic optimization belongs to the core of decision-making under uncertainty with solution methods ranging from theoretical techniques onto numerical approaches leveraging dynamic control and approximation of scenario trees (Shapiro, Dentcheva and Ruszczyński [51]). Clearly, the analytical solution can be found explicitly only in very exceptional cases due to variational form of such problems. More complex formulations can be handled by numerical solution techniques leading to the necessity to approximate the uncertainty set.

Dynamic optimization aims to find an approximate solution efficiently but usually neglects an important part of interstage dependencies while predicting future outcomes (Bellman [5], Bertsekas [7] and Dreyfus [17]). Methods on scenario trees can represent strong interdependencies between time stages being, however, subject to the curse of dimensionality (Fort and Pagés [21], Pflug and Römisch [41], Römisch [50]). To ease this problem, it is common to approximate the uncertainty by a tree with a smaller number of scenarios but the quality of the resulting decisions suffers significantly if Monte-Carlo simulations are applied for the approximation. This occurs because the Law of Large Numbers must be satisfied to guarantee the quality of decisions based on the Monte-Carlo sampling.

In this article, we consider stochastic processes given by *continuous-state* Fréchet probability distributions, estimated data-based and changing over time conditionally on new realizations (Mirkov and Pflug [36], Mirkov [37]). Based on these estimates, we approximate stochastic processes by finitely-valued *scenario trees* (Heitsch and Römisch [28], Pflug and Pichler [43]), which we directly use to solve multi-stage stochastic optimization problems. In order to assess the quality of such approximation techniques in multi-stage stochastic optimization problems, we formulate two optimization problems: The *initial* multi-stage stochastic optimization program, which is formulated in a continuous form, and the *approximate* problem, which is finite and discrete (Pflug and Römisch [41], Pflug and Pichler [43], Römisch [50]). The distance between these problems determines the approximation error. As the problems are generally defined on different probability spaces, the approximation error can bounded via the use of the *nested distance* introduced by Pflug and Pichler [43, 44]. The distance generalizes the well-known Kantorovich-Wasserstein distance from a single- to a multi-stage case (see Kantorovich [30], Pflug and Picher [43], Villani [59]). Various contorted nested distances, e.g., the weighted measure accounting for heavy-tail distributions, are defined for special situations as in Birghila and Pflug [11].

Importantly, the minimization of the nested distance over a finite nested distribution with a given tree structure would lead to the *optimal scenario tree* able to provide fine approximation of the initial optimization problem. Nevertheless, the minimization of the nested distance between a continuous-state stochastic process and a finitely valued scenario tree is an intractable optimization problem. Instead, current research focuses on the introduction of various bounds allowing for efficiency (see Deng et al. [16], Pflug and Pichler [43, 44]). In particular, the sum of Kantorovich-Wasserstein distances between stage-wise probability measures provides one of the upper bounds on the nested distance. However, each stage-wise minimization of the Kantorovich-Wasserstein distance is at least #P-hard (see Taşkesen [55]) and the optimal value is strongly dependent on the scenario tree structure, while slight changes lead to an increase or to a decrease of the nested distance and its upper bounds. In our article, we first of all demonstrate that both the nested distance and the approximation error converge to zero inversely proportional to the minimal number of quantizers, $n$, among all conditional distribution at a given time stage. Afterwards, we enhance the stage-wise optimal scenario generation for Fréchet distributions. In particular, we propose an extremely efficient method based on changing medians to approximate optimal quantizers at later stages using the optimal quantizers at the initial stage only. This allows to avoid sequential solution of the optimal transport problems in multi-stage optimization.

Next, we combine the stage-wise method for optimal scenario quantization with dynamic programming as of Bellman [5], Bertsekas [7] and Dreyfus [17], who expressed the optimal policy in terms of an optimization problem with iteratively evolving value function (the optimal *cost-to-go* function). We build on recent works of Ermoliev et al. [20], Bertsekas [8], Keshavarz [31], Hanasusanto and Kuhn [25], Powell [46], Li, [34], Sun [53], who demonstrate that the evaluation of optimal cost-to-go functions, involving multivariate conditional expectations, is a computationally complex procedure. This leads to the necessity to develop accurate and numerically efficient algorithms for multi-stage stochastic optimization. For example, the article of Hanasusanto and Kuhn [25] proposes an efficient piece-wise linear approximation of value functions, which can be incorporated directly into the dynamic programming. However, such an approximation does not necessarily guarantee the achievement of the global solution in the initial problem. Moreover, according to Baillon et al. [3], misperceived probabilities of considered scenarios often lead to a bias in optimal strategies. In our article, we use optimal weights for the value function approximation at each stage of the scenario tree. We prove and implement the conditions on the values function approximation able to guarantee the global solution. In particular, we demonstrate that monotonicity and convexity/concavity properties of value functions stay recursive in the dynamic programming. Thus, we introduce semidefinite constraints able to guarantee the global solution of the initial optimization problem. By this, the value function approximation becomes an SDP (semidefinite programming) problem.

Note that one needs to resolve two problems in order to propose an accurate and efficient algorithm for dynamic solution of multi-stage problems using optimal quantizers. Firstly, due to the convergence of the nested distance in the number of quantizers $n$ of each conditional distribution, one needs to be able to construct and store scenario trees of large sizes very efficiently. Secondly, one needs to reduce the number of problems to be solved in the dynamic programming as much as possible to preserve accuracy and ease the curse of dimensionality. The first problem is partly resolved for the Gaussian distribution due to the fact that 1) its conditional form is known explicitly (see Lipster and Shiryayev [35]), 2) the quantizers can be approximated linearly using the quantizers of the standard normal distribution, the mean and the variance and 3) there are accurate methods that estimate the multivariate probability integrals (Szántai [54]). Furthermore, the lognormal distribution can be approximated using the exponents of the Gaussian quantizers (see Timonina [56]), which makes this case simplistic as well. Clearly, a general distribution can be approximated by a finite Gaussian mixture in a semi-parametric way (Fruhwirth-Schnatter et al. [22]). However, an additional approximation solved for each subtree does not make the problem more efficient. Thus, in this article, we propose a quick way to approximate optimal quantizers for conditional Fréchet distributions, which are often used for flood risk modeling (see Wemelsfelder [61]). We base our approximation on the famous work of Gumbel [24], which allows us to subdivide optimal quantizers into two groups dependent on the risk threshold.

We enhance efficiency of the dynamic programming via differentiating between methods used for higher- and lower-risk quantizers: for quantizers whose percentile is lower than the Probability of No Loss (PNL) (see Timonina et al. [57]), we use our fast estimation; for rare but damaging events, we use the quick estimation of full parameter set proposed by Gumbel [24]. Such a combination allows to enhance efficiency for up to $100\% \cdot \text{PNL}$ of value function approximations, where $\text{PNL} = 0.6779$ for the Austrian case study. Furthermore, for large sample sizes, $N$, the higher is the risk threshold, the lower is the time complexity of the dynamic programming: we demonstrate that the number of dynamic optimization problems becomes linear in the number of quantizers, $n$, which is a significant enhancement compared to at least $n^t$ dynamic problems solved at stage $t$ without our approximation.

Overall, we demonstrate that our stage-wise quantization method combines high approximation quality with computational efficiency. It leads to the optimal weighting of value function estimates in the backward step and enhances both accuracy and efficiency of the decision-making due to the use of smaller-size but accurate scenario trees. This is due to the fact that optimal methods allow for using significantly smaller-size scenario trees in comparison to the necessity of sampling a large number of points in case of Monte-Carlo simulations. Our method accounts for the information about the past due to the use of conditional distribution functions at later stages of scenario trees. To improve the quantization even further, one could also incorporate information about feasible future scenarios proposing a combination of the forward procedure with the backward step on the tree.

A wide variety of applications of multi-stage stochastic optimization includes such fields as financial and investment planning, inventory control, and also natural disaster risk-management. We apply the developed algorithm to multi-period and multi-region risk-management of flood events in Austria, constructing a governmental decision-making problem, which takes the impact of natural hazards into account. For this, we use marginal losses which are simulated using the LISFLOOD hydrological model and an economic damage model as in the works of Van der Knijff, Younis and De Roo [32], as well as Rojas et al. [49]. We couple marginal losses in localities in such a way that the large-scale Austrian probability distribution is not underestimated and fits the multi-regional data on losses. To assess the interregional dependency, we use the structured coupling approach as of Timonina et al. [57]. We formulate the optimization in terms of a multi-stage stochastic problem with disaster losses described by a data-driven Fréchet distribution, which is, however, subject to change dependent on new event realizations. As a result, we find optimal strategies using the developed solution method avoiding the construction of a full scenario tree.

The article proceeds as follows: Section 2 describes optimal scenario generation method and the quantizer interpolation procedure. In Section 3 we rewrite the initial multi-stage stochastic optimization problem in the dynamic programming form and combine optimal scenario generation with backtracking solution methods of multi-stage stochastic optimization. In Section 4 we apply the algorithms to the problem of risk-management of flood events in Austria.

## 2. Mathematical framework

Given a discrete time horizon $T$ (i.e., $t \in \{1,...,T\}$), we consider a multi-stage stochastic optimization program in the following form with loss/profit function $H(x,\xi) = h_0(x_0) + \sum_{t=1}^{T} h_t(x^t,\xi^t)$ (e.g., Pflug and Römisch [41], Pflug [42], Pflug and Pichler [43, 44]):

$$\inf_{x \in \mathbb{X}, x \lhd \mathscr{F}} \mathbb{E}\left[ H(x,\xi) = h_0(x_0) + \sum_{t=1}^{T} h_t(x^t,\xi^t) \right]. \tag{1}$$

Here, $\xi = (\xi_1,...,\xi_T)$ is a $r_1$-dimensional continuous-state stochastic process ($\xi_t \in \mathbb{R}^{r_1}$, $\forall t = 1,...,T$) defined on the probability space $(\Omega,\mathscr{F},P)$ and $\xi^t = (\xi_1,...,\xi_t)$ is its history up to time $t$. Each random variable $\xi_t$ is measurable with respect to $\sigma$-algebra $\mathscr{F}_t$, $\forall t = 1,...,T$, where $\mathscr{F} = (\mathscr{F}_1,...,\mathscr{F}_T)$ is a filtration on the space $(\Omega,\mathscr{F},P)$. We denote this measurability as $\xi \lhd \mathscr{F}$ and add the trivial $\sigma$-algebra $\mathscr{F}_0 = \{\emptyset,\Omega\}$ as the first element of the filtration $\mathscr{F}$. Importantly, the sequence of decisions $x = (x_0,...,x_T)$ ($x_t \in \mathbb{R}^{r_2}$, $\forall t = 0,...,T$) with history $x^t = (x_0,...,x_t)$ must also satisfy the constraint $x \lhd \mathscr{F}$, which is referred to as *the non-anticipativity condition* (e.g., Pflug [42], Pflug and Pichler [43, 44]). Indeed, this means that only those decisions are feasible, which are based on the information available at the particular time. Additionally, $\mathbb{X}$ denotes the set of constraints on $x$ other than the non-anticipativity constraints. The problem (1) is stated for a general continuous-state stochastic process $\xi$. In this article, however, we focus on the case, when each random variable $\xi_t$ follows a Fréchet distribution $P_t(\xi_t|\xi^{t-1}) = \exp\left[ -\left(\frac{u_t-\varepsilon_t}{\xi_t-\varepsilon_t}\right)^{\frac{1}{\lambda_t}} \right]$, $\lambda_t > 0$, $\varepsilon_t \leq \xi_t < \infty$, where $u_t$ and $\varepsilon_t$ are location parameters and $\lambda_t$ is a shape parameter changing in time dependent on the observations of the stochastic process until time $t$.

The variational form of the optimization problem (1) does not allow to solve it in a closed form, which implies the necessity to approximate the problem for a numerical solution. The approximated problem (2) can be written in the form with the loss/profit function $H(\widetilde{x},\widetilde{\xi}) = h_0(\widetilde{x}_0) + \sum_{t=1}^{T} h_t(\widetilde{x}^t,\widetilde{\xi}^t)$:

$$\inf_{\widetilde{x} \in \mathbb{X}, \widetilde{x} \lhd \widetilde{\mathscr{F}}} \mathbb{E}\left[ H(\widetilde{x},\widetilde{\xi}) = h_0(\widetilde{x}_0) + \sum_{t=1}^{T} h_t(\widetilde{x}^t,\widetilde{\xi}^t) \right], \tag{2}$$

where the stochastic process $\xi$ is replaced by a scenario process $\widetilde{\xi} = (\widetilde{\xi}_1,...,\widetilde{\xi}_T)$, such that $\widetilde{\xi}_t \in \mathbb{R}^{r_1}$, $\forall t = 1,...,T$ and $\widetilde{\xi}_t$ is discrete (i.e., $\widetilde{\xi}_t$ takes finite number of values $n_t$, $\forall t = 1,...,T$). Scenario process $\widetilde{\xi} = (\widetilde{\xi}_1,...,\widetilde{\xi}_T)$ can be defined on a probability space $(\widetilde{\Omega},\widetilde{\mathscr{F}},\widetilde{P})$, which is different from $(\Omega,\mathscr{F},P)$ (see Pflug and Römisch [41], Pflug [42], Pflug and Pichler [43, 44]).

The distance between problems (1) and (2) determines the approximation error. Until recently, the distance between the initial problem (1) and its approximation (2) was defined only if both processes $\xi$ and $\widetilde{\xi}$ and both filtrations $\mathscr{F}$ and $\widetilde{\mathscr{F}}$ were defined on the same probability space $(\Omega,\mathscr{F},P)$, meaning that the approximation error was measured as a filtration distance. The introduction of the concept of the *nested* *distribution* (see Pflug [42], Pflug and Pichler [43, 44]), containing in one mathematical object the scenario values as well as the structural information under which decisions have to be made, allowed to bring the problem to the purely distributional setup. The *nested distance* between these distributions was first introduced by Pflug and Pichler [43, 44] and turned out to be a multi-stage generalization of the well-known *Kantorovich-Wasserstein distance* defined for single-stage problems (see Kantorovich [30], Pflug and Pichler [43, 44], Villani [59]). Minimizing the nested distance, one can enhance the quality of the approximation and, hence, the solution accuracy. Numerical methods for the nested distance reduction are based on the construction of tight upper bounds and their minimization (e.g., Fort and Pagés [21], Pflug and Römisch [41], Römisch [50]). This includes the stage-wise minimization of the Kantorovich-Wasserstein distance between measures sitting at each stage of the scenario tree (e.g., Villani [59], Timonina [56]).

In this article, we combine the stage-wise method for distributional quantization of Fréchet distributions with a backtracking solution algorithm, which is based on the dynamic programming principle (e.g., Ermoliev, Marti and Pflug [20], Hanasusanto and Kuhn [25]). This leads to the optimal weighting of value function estimates in the backward step and enhances both accuracy and efficiency of the decision-making due to the use of smaller-size but accurate scenario trees. Furthermore, via the introduction of some particular constraints in the value function approximation, we are able to guarantee the global solution of the problem under convexity and monotonicity conditions. We introduce the general method for optimal scenario tree quantization in Section 2.1. Grouping the quantizers of Fréchet distributions in Section 2.2, we efficiently combine them with the dynamic programming in Section 3.

### 2.1. Stage-wise optimal scenario tree approximation

A wide range of forward-looking numerical methods for the solution of multi-stage stochastic optimization problems is based on the approximation of stochastic process $\xi = (\xi_1, ..., \xi_T)$ by *finitely valued* scenario trees (see Appendix A.1). To achieve the best approximation quality using Monte-Carlo sampling, one needs the Law of Large Numbers to hold. Nevertheless, instead of using Monte-Carlo samples, it is very natural to minimize the distance between the continuous distribution $P_t(\xi_t|\xi^{t-1})$ and a discrete probability measure $\widetilde{P}_t(\widetilde{\xi}_t|\widetilde{\xi}^{t-1}) = \sum_{i=1}^n \widetilde{p}_t^i(\widetilde{\xi}^{t-1})\delta_{z_t^i(\widetilde{\xi}^{t-1})}$ sitting on $n$ points $z_t^i(\widetilde{\xi}^{t-1})$, $\forall i = 1, ..., n$, with corresponding probabilities $\widetilde{p}_t^i(\widetilde{\xi}^{t-1}), \forall i$. In the following, we consider the well-known *Kantorovich-Wasserstein distance* between probability measures (Kantorovich [30], Villani [59]) and the connection between its minimization and the minimization of the approximation error between problems (1) and (2).

The Kantorovich-Wasserstein distance is used in a wide variety of areas, allowing for accurate modeling of scenarios and ambiguity sets in distributionally robust optimization problems (e.g., Blanchet [12]).

DEFINITION 1. The Kantorovich distance between probability measures $P$ and $\widetilde{P}$ is defined as

$$d_{KA}(P, \widetilde{P}) = \inf_{\pi} \left\{ \int_{\Omega \times \widetilde{\Omega}} d(w, \widetilde{w}) \pi[dw, d\widetilde{w}] \right\}, \tag{3}$$
$$\text{subject to } \pi[\cdot \times \widetilde{\Omega}] = P(\cdot) \text{ and } \pi[\Omega \times \cdot] = \widetilde{P}(\cdot),$$

where $d(w,\widetilde{w})$ is the cost function for the transportation of $w \in \Omega$ to $\widetilde{w} \in \widetilde{\Omega}$.

Importantly, the Kantorovich-Wasserstein distance accounts for available stage-wise information only, though more information is available on a tree structure. Therefore, one cannot guarantee that the stage-wise minimization of the Kantorovich-Wasserstein distance always results in the minimal approximation error between problems (1) and (2). To overcome this dilemma, Pflug and Pichler introduced the concept of *nested distributions* $\mathbb{P} \sim (\Omega, \mathscr{F}, P, \xi)$ and $\widetilde{\mathbb{P}} \sim (\widetilde{\Omega}, \widetilde{\mathscr{F}}, \widetilde{P}, \widetilde{\xi})$ (see the article [43]), which contain information about both processes $\xi$ and $\widetilde{\xi}$ and all underlying interdependencies. The *nested (multi-stage) distance* between nested distributions is defined in a purely distributional setup (Pflug and Pichler [44]). In our article, the nested distance is denoted by $dl(\mathbb{P},\widetilde{\mathbb{P}})$, where $\mathbb{P}$ refers to the continuous nested distribution of the initial problem (1) and $\widetilde{\mathbb{P}}$ corresponds to the discrete nested distribution, which is the scenario tree approximation of the problem (2).

DEFINITION 2. The nested (multi-stage) distance (see Pflug and Pichler [43, 44]) is defined as

$$\begin{aligned} dl(\mathbb{P},\widetilde{\mathbb{P}}) = \ & \inf_{\pi} \left( \int d(w,\widetilde{w}) \pi(dw, d\widetilde{w}) \right), \\ & \text{subject to } \mathbb{P} \sim (\Omega, \mathscr{F}, P, \xi),\ \widetilde{\mathbb{P}} \sim (\widetilde{\Omega}, \widetilde{\mathscr{F}}, \widetilde{P}, \widetilde{\xi}) \\ & \pi[A \times \widetilde{\Omega} | \mathscr{F}_t \otimes \widetilde{\mathscr{F}_t}](w,\widetilde{w}) = P(A|\mathscr{F}_t)(w),\ (A \in \mathscr{F}_T, 1 \leq t \leq T), \\ & \pi[\Omega \times B | \mathscr{F}_t \otimes \widetilde{\mathscr{F}_t}](w,\widetilde{w}) = \widetilde{P}(B|\widetilde{\mathscr{F}_t})(\widetilde{w}),\ (B \in \widetilde{\mathscr{F}_T}, 1 \leq t \leq T). \end{aligned} \tag{4}$$

Under the assumption of Lipschitz-continuity of the loss/profit function $H(x,\xi)$ with the Lipschitz constant $L_1$, the nested distance (4) establishes an upper bound on the approximation error, which implies that $|v(\mathbb{P}) - v(\widetilde{\mathbb{P}})| \leq L_1 dl(\mathbb{P},\widetilde{\mathbb{P}})$ (Pflug and Pichler [44]). Here, value functions $v(\mathbb{P})$ and $v(\widetilde{\mathbb{P}})$ correspond to optimal solutions of the multi-stage problems (1) and (2). If the Lipschitz property also holds for the stage-wise Kantorovich-Wasserstein distances with constants $L_t$ $\forall t = 2,...,T$, i.e., $d_{KA}(P_t(\cdot|w^{t-1}), P_t(\cdot|\widetilde{w}^{t-1})) \leq L_t d(w^{t-1}, \widetilde{w}^{t-1})$, where $(w_1,...,w_T)$ and $(\widetilde{w}_1,...,\widetilde{w}_T)$ are points from the sample space with histories $w^{t-1} = (w_1,...,w_t)$ and $\widetilde{w}^t = (\widetilde{w}_1,...,\widetilde{w}_t)$ correspondingly, the following can be derived based on Lemma 3.11 in the work of Pflug and Pichler [43] with $d_{KA}(P_t,\widetilde{P}_t) = \sup_{w^{t-1}} d_{KA}(P_t(\cdot|w^{t-1}), \widetilde{P}_t(\cdot|w^{t-1}))$:

$$dl(\mathbb{P},\widetilde{\mathbb{P}}) \leq \sum_{t=1}^{T} d_{KA}(P_t,\widetilde{P}_t) \prod_{s=t+1}^{T} (L_s + 1). \tag{5}$$

Now, we let $n_t$, $\forall t = 1,...,T$ be the total number of scenarios at the stage $t$ and $n_t^i$ ($\forall i = 1,...,n_{t-1}, \forall t = 2,...,T$) be the number of quantizers corresponding to the $n_{t-1}$ conditional distributions sitting at the stage $t$ ($\forall t = 2,...,T$). We note that $n_t = \sum_{i=1}^{n_{t-1}} n_t^i, \forall t > 1$. According to the assymptotic quantization error theorem of Graf and Luschgy [23], there exists $\widetilde{P}_t^{n_t^i}$ ($\forall t = 1,...,T$) sitting on $n_t^i$ discrete points, such that

$$d_{KA}(P_t,\widetilde{P}_t^{n_t^i}) \leq c(n_t^i)^{-\frac{1}{r_1}} \leq c \left( \min_{i \in \{1,..,n_{t-1}\}} n_t^i \right)^{-\frac{1}{r_1}} \leq c \left( \min_{\substack{i \in \{1,..,n_{t-1}\} \\ t \in 1,...,T}} n_t^i \right)^{-\frac{1}{r_1}}. \tag{6}$$

Thus, the nested distance convergence result follows:

$$|v(\mathbb{P}) - v(\widetilde{\mathbb{P}})| \leq cL_1 n^{-\frac{1}{r_1}} \sum_{t=1}^{T} \prod_{s=t+1}^{T} (L_s + 1), \tag{7}$$

where $n = \min_{i \in \{1,..,n_{t-1}\}, t \in \{1,...,T\}} n_t^i$.

Independent of the tree structure, minimizing the Kantorovich-Wasserstein distance at each stage of the scenario tree by finding $n = \min_{i \in \{1,..,n_{t-1}\}, t \in \{1,...,T\}} n_t^i$ optimal supporting points and their probabilities leads to the optimal quantization $\widetilde{P}_t^n$ which satisfies (6) and guarantees a speed of convergence at least as good as in the upper bound (7). This allows to enhance accuracy of the solution of multi-stage stochastic optimization problems with respect to the well-known Monte-Carlo (random) scenario generation. Note that the optimal quantization finds *n optimal supporting points* $\xi_t^{*i},\ i = 1,...,n$ of conditional distribution $P_t(\xi_t|\xi^{t-1}),\ \forall t = 2,...,T$ by minimizing the following distance in $z_t^i,\ \forall i$:

$$\mathscr{D}\left(z_t^1,...,z_t^n\right) = \int \min_i d\left(w_t, z_t^i\right) P_t(dw_t|w^{t-1}), \tag{8}$$

where $d\left(w_t, z_t^i\right)$ is the Euclidean distance between points $w_t$ and $z_t^i$. At stage $t = 1$, optimal quantization is based on the unconditional distribution $P_1(\xi_1)$. Given the locations of the supporting points $\xi_t^{*i}$, their probabilities $p_t^{*i}$ are calculated by the minimization of the Kantorovich-Wasserstein distance between the measure $P_t(\xi_t|\xi^{t-1})$ and its discrete approximation:

$$\min_{p_t^{*i},\,\forall i} d_{KA}\left(P_t(\xi_t|\xi^{t-1}), \sum_{i=1}^{n} p_t^{*i} \delta_{\xi_t^{*i}}\right), \tag{9}$$

where the hardness of such semi-discrete optimal transport problems is studied by Taşkesen et al. [55] in Theorem 2.2. Further in our article, we propose a quick way to approximate optimal quantizers for conditional Fréchet distributions, which are often used for flood risk modeling. We proceed with the approximation of optimal supporting points for the Fréchet distribution in Section 2.2.

### 2.2. Parameters of Fréchet distributions and quantizer grouping

In the famous work of Gumbel [24], the author proposes a quick method to estimate parameters of the Fréchet distribution using the median $\hat{x}_s$, the smallest value $\hat{x}_1$ and the largest value $\hat{x}_N$ in a sample of size $N \geq 3$ of a random variable $\xi$. In particular, they consider the distribution function $F(\xi) = \exp\left[-\left(\frac{u-\varepsilon}{\xi-\varepsilon}\right)^{\frac{1}{\lambda}}\right]$, $\lambda > 0$, $\varepsilon \leq \xi < \infty$, where $u$ and $\varepsilon$ are location parameters and $\lambda$ is a shape parameter.

The following estimates are derived by Gumbel [24] for the lower limit $\varepsilon$ and the difference $u - \varepsilon$:

$$\hat{\varepsilon} = \frac{\hat{x}_s N^{\lambda} - \hat{x}_N}{N^{\lambda} - 1},, \quad \hat{u} - \hat{\varepsilon} = (\hat{x}_s - \hat{\varepsilon})(\log 2)^{\lambda}. \tag{10}$$

To estimate the shape parameter $\lambda$, one needs to solve the following implicit equation with respect to $\hat{\lambda}$:

$$\frac{\hat{x}_N - \hat{x}_s}{\hat{x}_s - \hat{x}_1} = \frac{N^{\hat{\lambda}} - 1}{1 - f(\hat{\lambda}, N)(\log 2)^{\hat{\lambda}}} = g(\lambda, N), \tag{11}$$

where $f(\lambda,N) = \left\{ -\log[1-(0.5)^{\frac{1}{N}}] \right\}^{-\lambda}$ and where we denote the right-hand side by $g(\lambda,N)$. Figure 1 demonstrates the numerical estimate of the residual between the left and right sides of the equality (11) for different values of $\lambda$. The value of $\lambda$ with the minimal residual is the estimate $\hat{\lambda}$.

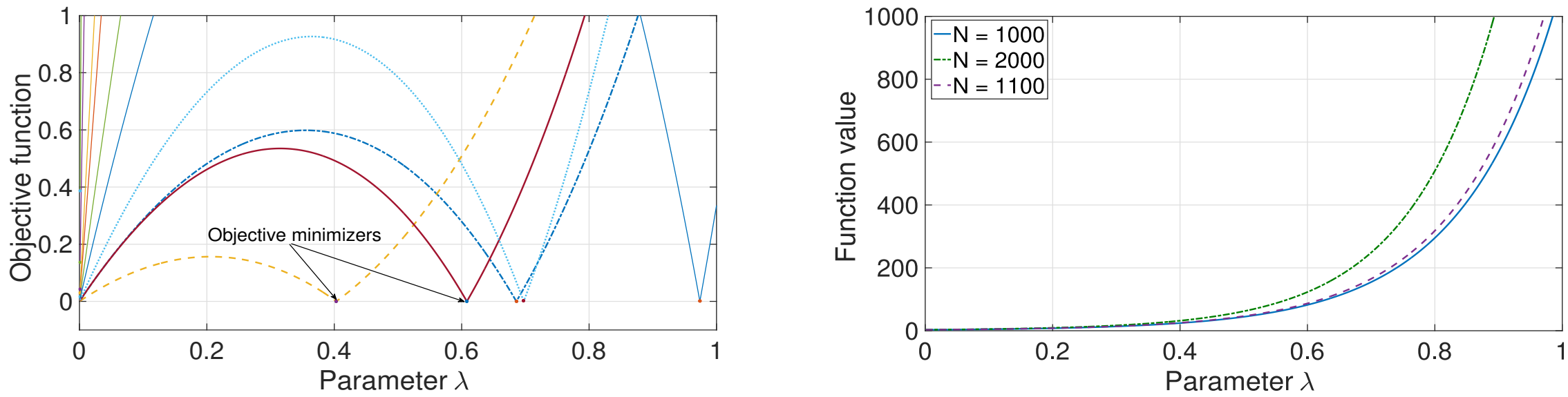

a. Minimization of the residual between the left and right sides of the equality (11).

b. Function value of the right-hand side in equality (11).

**Figure 1 Minimization of the residual between the left and right sides of the equality (11).**

As one can see in Figure 1a, the residual function is neither convex-concave nor smooth in $\lambda$. If one could assume that the shape parameter $\lambda$ stays constant over time, one could propose a quick method to update the parameters of the Fréchet distribution based on a new sample point. Indeed, given a $(N+1)$-th sample point, one would update the parameters in the following way:

$$\hat{\varepsilon}^{\text{new}} = \frac{\hat{x}_s^{\text{new}}(N+1)^{\lambda} - \hat{x}_{N+1}^{\text{new}}}{(N+1)^{\lambda}-1} = \frac{N^{\lambda}\hat{x}_s\left(\frac{\hat{x}_s^{\text{new}}}{\hat{x}_s}\right)(\frac{1}{N}+1)^{\lambda} - \hat{x}_N\left(\frac{\hat{x}_{N+1}^{\text{new}}}{\hat{x}_N}\right)}{N^{\lambda}(\frac{1}{N}+1)^{\lambda}-1} \approx \frac{N^{\lambda}\hat{x}_s\left(\frac{\hat{x}_s^{\text{new}}}{\hat{x}_s}\right) - \hat{x}_N\left(\frac{\hat{x}_{N+1}^{\text{new}}}{\hat{x}_N}\right)}{N^{\lambda}-1}, \quad (12)$$

$$\frac{\hat{u}^{\text{new}}}{\hat{u}} = \frac{\hat{\varepsilon}^{\text{new}} + (\log 2)^{\lambda}(\hat{x}_s^{\text{new}} - \hat{\varepsilon}^{\text{new}})}{\hat{\varepsilon} + (\log 2)^{\lambda}(\hat{x}_s - \hat{\varepsilon})} = \frac{\hat{\varepsilon}^{\text{new}}}{\hat{\varepsilon}} \cdot \frac{1+(\log 2)^{\lambda}\left(\frac{\hat{x}_s^{\text{new}}}{\hat{\varepsilon}_s^{\text{new}}}-1\right)}{1+(\log 2)^{\lambda}\left(\frac{\hat{x}_s}{\hat{\varepsilon}_s}-1\right)}, \quad (13)$$

where $\hat{\varepsilon}^{\text{new}}$ and $u^{\widehat{\text{new}}}$ are the new estimates of location parameters (note that the superscript *new* implies the updated value). Importantly, the shape parameter $\lambda$ has low sensitivity to large changes in the function $g(\lambda,N)$ if the value of $\lambda$ is high (see Figure 1b). Furthermore, both sides of the equality (11) must follow the same trend if the new sample point occurs: if the left-hand side of the equality grows, the right-hand side must grow as well; oppositely, when the left-hand side value drops, the right-hand side must also decrease so that the equality is satisfied. Thus, consider the following three cases for a new sample point:

***Case 1 - Smaller losses:*** As the median $\hat{x}_s$ is a 50%-percentile of a distribution function, there is a 50% chance that the new sample point falls below $\hat{x}_s$. In this case, $\frac{\hat{x}_s^{\text{new}}}{\hat{x}_s} \leq 1$, while $\frac{\hat{x}_{N+1}^{\text{new}}}{\hat{x}_N} = 1$ as the largest value is not influenced. Based on the equation (12), we can state:

$$\hat{\varepsilon}^{\text{new}} \leq \hat{\varepsilon}\frac{\hat{x}_s^{\text{new}}}{\hat{x}_s} \leq \hat{\varepsilon}.$$

Next, as $\hat{x}_1$ (which is equal to zero for flood losses) and $\hat{x}_N$ are not changing, the left-hand side of equality (11) increases. As the function $g(\lambda,N)$ has low sensitivity to small changes in the number of sample points for large $N$ (note that its elasticity in $N$ is proportional to $\frac{\lambda}{N^{\lambda}-1}$), the right-hand side would only increase if $\lambda$ is constant or increasing (see Figure 1b). Assuming no change in the shape parameter, we select the maximal value of $\hat{\varepsilon}^{\text{new}}$ in order to avoid underestimation of extremes, which implies the updates $\hat{\varepsilon}^{\text{new}} = \hat{\varepsilon}\frac{\hat{x}_s^{\text{new}}}{\hat{x}_s}$ with the constant $\lambda$.

***Case 2 - Significant losses:*** There is also a 50% chance that the new sample point falls above $\hat{x}_s$. If the largest value in the sample is not influenced, one can write $\frac{\hat{x}_s^{\text{new}}}{\hat{x}_s} \geq 1$ and $\frac{\hat{x}_{N+1}^{\text{new}}}{\hat{x}_N} = 1$. It yields:

$$\hat{\varepsilon} \leq \hat{\varepsilon}\frac{\hat{x}_s^{\text{new}}}{\hat{x}_s} \leq \hat{\varepsilon}^{\text{new}}.$$

The left-hand side of equality (11) decreases, while the right-hand side, given that $\lambda$ is constant, increases (see Figure 1b). The equality (11) is clearly not satisfied under the condition of constant $\lambda$ in this case.

***Case 3 - The largest loss:*** There is a small chance that the new sample point falls above $\hat{x}_N$. If this happens, the largest observation in the sample is updated, i.e., $\frac{\hat{x}_s^{\text{new}}}{\hat{x}_s} \geq 1$ and $\frac{\hat{x}_{N+1}^{\text{new}}}{\hat{x}_N} \geq 1$. The left-hand side of the equality (11) can increase or decrease, dependent on the behavior of the nominator. Clearly, in this case the assumption of a constant parameter $\lambda$ is not valid any more as the event is large enough to influence it.

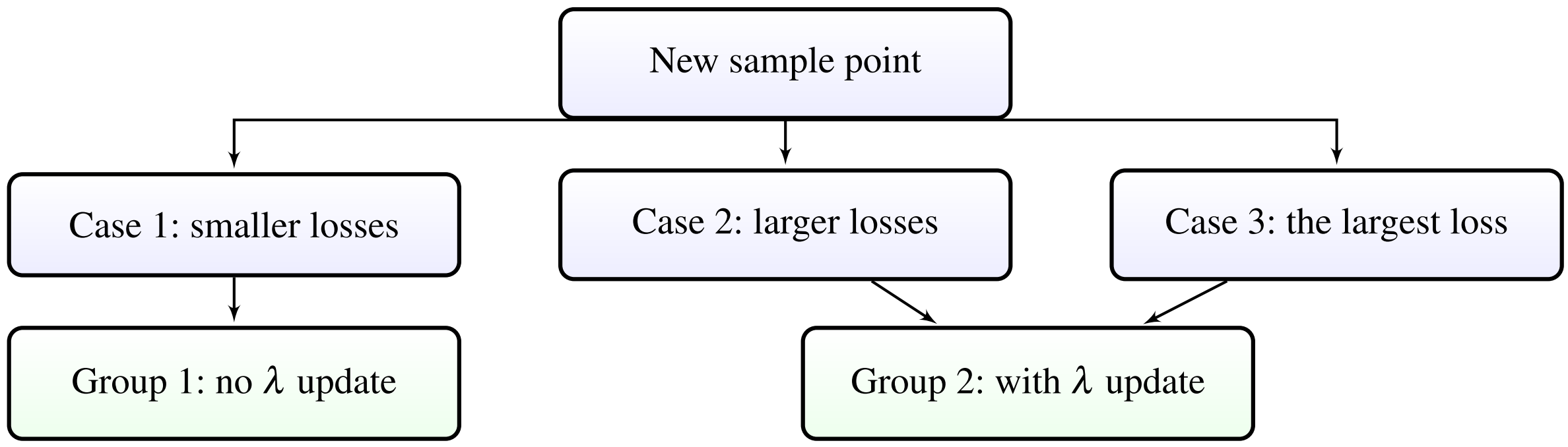

**Figure 2** **Grouping of new samples based on the estimated median.**

Overall, based on these cases, we subdivide samples into two groups shown in Figure 2. As the shape parameter $\lambda$ changes in case of larger flood events, the decision-maker is better off using the complete procedure (10-11) for Group 2. However, the updates (14) can be used for the case of losses in Group 1:

$$\hat{\varepsilon}^{\text{new}} = \hat{\varepsilon}\frac{\hat{x}_s^{\text{new}}}{\hat{x}_s},\ \hat{u}^{\text{new}} = \hat{u}\frac{\hat{\varepsilon}^{\text{new}}}{\hat{\varepsilon}}. \tag{14}$$

Note that cases 1)-3) use the percentile of 50%, which corresponds to the probability of 0.5 (referred to as the *risk threshold*) that a new sample point falls below the median. However, a larger risk threshold can be used for modeling of flood losses in Austria. Indeed, based on the raw data on natural disasters and their times, one could estimate the maximal loss $L$, whose appearance would not influence the shape of the distribution. This could be done by appending different values of $L$ and testing the hypothesis that the shape parameter stays constant given a new entry. The largest $L$ for which the hypothesis holds true can serve to compute the risk threshold $F(L)$.

Differently, one could find the risk threshold based on the Probability of No Loss defined as PNL $= 1 - P(N_L > 1)$, where $N_L$ is the number of *exceedance events*, i.e., events causing damages and economic losses (see Calenda [13] and Timonina et al. [57]). Employing the estimation method of Timonina et al. [57] for computing the PNL, we conclude that the PNL = 0.6779 for Austria, which is higher than 0.5 (see Appendix A.2). Thus, if the new sample point falls below PNL in probability, we still use the simplified procedure $\hat{\varepsilon}^{\text{new}} = \hat{\varepsilon}\frac{\hat{x}_s^{\text{new}}}{\hat{x}_s}$ to update parameters of the loss distribution (Figure 3).

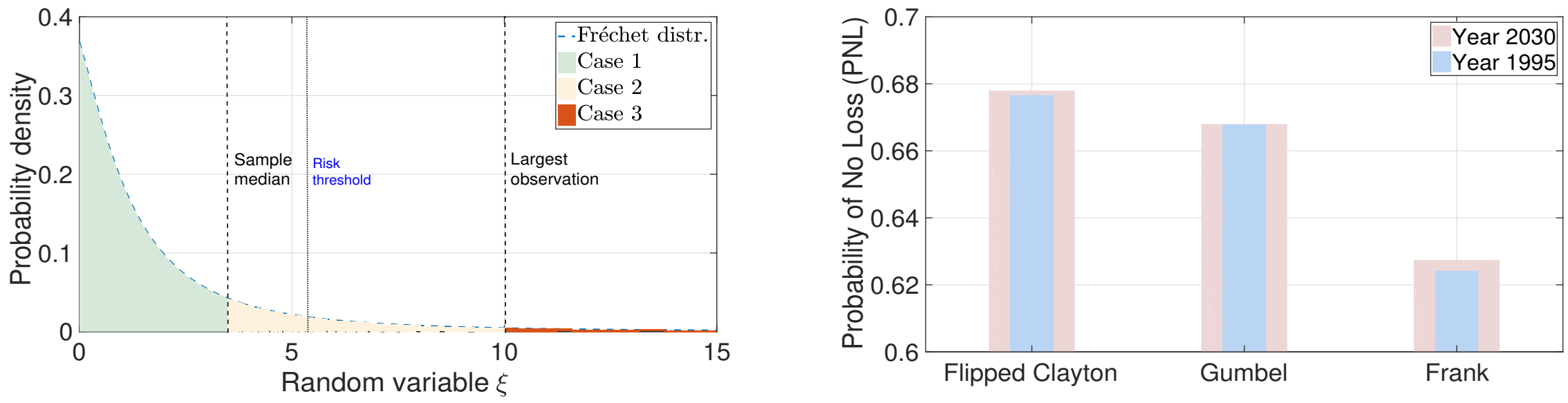

a. Grouping of samples. b. Probability of No Loss.

**Figure 3 Dividing Fréchet-distributed samples into groups based on PNL.**

The updates (14) allow us to preserve accuracy and enhance efficiency in up to 100·PNL% scenario quantization problems, while the actual gain depends on the number of distribution quantizers $n$ and on their probabilities. Furthermore, updates (14) allow to avoid solving multiple problems of the type (1) for generating optimal quantizers for random variables $\xi_t$ at each stage. Instead, one can use optimal quantizers for the random variable $\xi_{t-1}$ multiplied by $\frac{\hat{x}_s^{\text{new}}}{\hat{x}_s}$. Indeed, the distribution of the random variable $\frac{\hat{x}_s^{\text{new}}}{\hat{x}_s}\xi_{t-1}$ is Fréchet$(\lambda, \frac{\hat{x}_s^{\text{new}}}{\hat{x}_s}(u-\varepsilon), \frac{\hat{x}_s^{\text{new}}}{\hat{x}_s}\varepsilon)$, which is exactly the case with parameters (14). Thus, the following holds:

$$\mathscr{D}\left(\xi_t^{*1}, \ldots, \xi_t^{*n}\right) = \int \min_i d\left(\frac{\hat{x}_s^{\text{new}}}{\hat{x}_s} w_{t-1}, \frac{\hat{x}_s^{\text{new}}}{\hat{x}_s}\xi_{t-1}^{*i}\right) d\exp\left[-\left(\frac{w_t - \frac{\hat{x}_s^{\text{new}}}{\hat{x}_s}\varepsilon}{\frac{\hat{x}_s^{\text{new}}}{\hat{x}_s}(u-\varepsilon)}\right)^{-\frac{1}{\lambda}}\right] = \frac{\hat{x}_s^{\text{new}}}{\hat{x}_s}\mathscr{D}(\xi_{t-1}^{*1}, \ldots, \xi_{t-1}^{*n}),$$

while the optimal probabilities[1] do not change in line with the breakpoints discussed in Timonina [56].

Considering a path of optimal quantizers $(\xi_1^{*i}, \xi_2^{*i}, \ldots, \xi_t^{*i})$, where the index $i$ corresponds to the $i$-th quantizer of the distribution $P_t(\xi_t|\xi^{t-1})$, we note that if the $i$-th quantizer at stage $t$ belongs to Group 1, then the $i$-th conditional quantizer at stage $t+1$ also belongs to Group 1, because the optimal probabilities (9) $\sum_{j=1}^{i} p_t^{*j} \leq PNL$ do not change in time. This implies the recursion (15):

$$\xi_t^{*i} = \frac{\hat{x}_{s,t-1}^i}{\hat{x}_{s,t-2}^i}\xi_{t-1}^{*i} = \frac{\hat{x}_{s,t-1}^i}{\hat{x}_{s,t-3}^i}\xi_{t-2}^{*i} = \ldots = \frac{\hat{x}_{s,t-1}^i}{\hat{x}_{s,0}^i}\xi_1^{*i}, \text{ where } \xi_1^{*i} \leq \hat{x}_s \tag{15}$$

Therefore, it is enough to solve the problem (8) only once at stage $t = 1$, adapting the optimal supporting points linearly at all further stages for Group 1. This adaptation is one of the main benefits of our method for Fréchet distribution quantization and has not been studied in the literature before, though it provides fine and extremely efficient ground for dynamic programming for natural disaster risk management. Note that using the obtained quantizers and knowing the threshold between Groups 1 and 2 (or, Cases 1-3), one can easily differentiate between scenario tree nodes and, importantly, between methods used for the solution of the subproblems at these nodes.

[1] The optimal probabilities $p_t^{*i}$, $\forall i = 1, \ldots, n$ can be found as $p_t^{*i} = F(q_t^i) - F(q_t^{i-1})$ using the intervals $I_i = [q_t^{i-1}, q_t^i]$ with the breakpoints $q_t^0 = -\infty$, $q_t^n = \infty$ and $q_t^i = \frac{1}{2}\frac{\hat{x}_s^{\text{new}}}{\hat{x}_s}(\xi_{t-1}^{*i} + \xi_{t-1}^{*i+1})$ $\forall i = 1, \ldots, n-1$.

## 3. Dynamic programming with recursive quantizer updates

The idea of the dynamic programming method goes back to pioneering papers of Bellman [5], Bertsekas [7] and Dreyfus [17], who expressed the optimal policy in terms of an optimization problem with iteratively evolving *value function* (the optimal *cost-to-go* function). These foundational works provide us with the theoretical framework for rewriting time separable multi-stage stochastic optimization problems in the dynamic form. More recent works of Bertsekas [8], Keshavarz [31], Hanasusanto and Kuhn [25], Powell [46], Li, [34], Sun [53] are built on (i) the fact that the evaluation of optimal value functions, involving multivariate conditional expectations, is a computationally complex procedure and on (ii) the necessity to develop numerically efficient algorithms for multi-stage stochastic optimization. We follow this path and propose an accurate and efficient algorithm for the solution of multi-stage problems using optimal quantizers.

We consider *endogenous state variables* $s_t \in \mathbb{R}^{r_3}$ which accumulate all the decision-dependent information about the past (see Shapiro et al. [51] for model state equations for linear optimization). We assume that the dimension $r_3$ of the endogenous variable $s_t$ does not change in time and that the variable obeys the recursion $s_{t+1} = g_t(s_t, x_t, \xi_{t+1})\ \forall t$. The optimal values of these problems are denoted by $V_t(s_t, \xi_t)$ and are referred to as *value functions*, as they depend on the state $s_t$ and on the stochastic process realization $\xi_t$. Starting with the last stage, one minimizes the objective $h_t(s_t, x_t, \xi_t) + \mathbb{E}\left[V_{t+1}(s_{t+1}, \xi_{t+1}) \mid \xi_t\right]$ under a given set of constraints $\mathbb{X}_t$ and $s_{t+1} = g_t(s_t, x_t, \xi_{t+1})$ at each time period $t$.

Generating $K$ paths for the endogenous state variable $(\hat{s}_1^k, \dots, \hat{s}_T^k)\ \forall k = 1, \dots, K$, the estimate of $\hat{V}_t(\hat{s}_t^k, \xi_t^{*i})$ is obtained by solving the following problem $\forall i:\ \xi_1^{*i} \leq \hat{x}_{s,0}$:

$$\begin{aligned} \hat{V}_t(\hat{s}_t^k, \xi_t^{*i}) = \ & \min_{x_t} \left\{ h_t(\hat{s}_t^k, x_t, \xi_t^{*i}) + \sum_{j \in \mathscr{L}_{t+1}^i} \left[\hat{V}_{t+1}(s_{t+1}^j, \xi_{t+1}^{*j})\right] p_{t+1}^{*j} \right\}, \\ & \text{subject to } x_t \in \mathbb{X}_t,\ x_t \lhd \mathscr{F}_t,\ s_{t+1}^j = g_t(\hat{s}_t^k, x_t, \xi_{t+1}^{*j}),\ \forall j \in \mathscr{L}_{t+1}^i, \\ & \xi_{t+1}^{*j} = \frac{\hat{x}_{s,t}^i(\hat{s}_t^k)}{\hat{x}_{s,0}^i} \xi_1^{*j},\ \forall j \in \mathscr{L}_{t+1}^i \quad \textbf{(quantizer updates)}, \end{aligned} \tag{16}$$

where $\mathscr{L}_{t+1}^i$ is the set of node indices at stage $t+1$ emerging from the node with index $i$ at stage $t$ in the scenario tree. This set is necessary in order to define the chosen subtree and, therefore, preserve the information structure. Further, $\hat{x}_{s,0}^i$ is the median of the sample, which is independent of any decisions taken as it has already realized. Differently, $\hat{x}_{s,t}^i(\hat{s}_t^k)$ is the median of the sample including the scenario realization $i$ until the stage $t$. It might be dependent on prior decisions $x^{t-1}$ and, thus, on the state $\hat{s}_t^k$, especially if the decision-maker can influence the uncertainty via some measures (e.g., constructing large-enough dams for protection against floods).

Importantly, the optimization problem (16) requires the optimal value $\hat{V}_{t+1}(s^j_{t+1}, \xi^{*j}_{t+1})$ to be evaluated at the point $s^j_{t+1}$, which does not necessarily coincide with trajectories $\hat{s}^k_{t+1}$, $\forall k$. In the work of Hanasusanto and Kuhn [25], the piece-wise linear value function approximation is used for these purposes. In our article, we approximate the value of $\hat{V}_{t+1}(s_{t+1}, \xi^{*j}_{t+1})$ continuously in $s_{t+1}$ under the assumption about convexity and monotonicity of functions $h_t(s_t, x_t, \xi_t)$, $g_t(s_t, x_t, \xi_{t+1})$ and $V_{t+1}(s_{t+1}, \xi_{t+1})$. This allows to preserve efficiency and guarantees the global solution of the optimization problem. Indeed, if convexity and monotonicity conditions hold for functions $h_t(s_t, x_t, \xi_t)$, $g_t(s_t, x_t, \xi_{t+1})$ and $V_{t+1}(s_{t+1}, \xi_{t+1})$, one can claim that the function $V_t(s_t, \xi_t)$ is also convex and monotone. Moreover, these properties stay recursive $\forall t = T, ..., 0$ (Lemma 1).

LEMMA 1. Let $\xi_0$ and $\xi_1$ be two dependent random variables defined on some probability space $(\Omega, \mathscr{F}, P)$. Consider a function $h(s, x, \xi_0) : \mathbb{R}^r \to \mathbb{R}^1$, which is jointly convex in $(s, x)$ and monotonically increasing as $s_1 \geq s_2 \Rightarrow h(s_1, x, \xi_0) \geq h(s_2, x, \xi_0)$, $\forall x, \xi_0$. Let $g(s, x, \xi_1) : \mathbb{R}^r \to \mathbb{R}^r$ be componentwise convex in both $s$ and $x$ and componentwise increasing in $s$. Further, let $V_1(y, \xi_1) : \mathbb{R}^r \to \mathbb{R}^1$ be convex in $y$ and monotonically increasing as $y_1 \geq y_2 \Rightarrow V_1(y_1, \xi_1) \geq V_1(y_2, \xi_1)$, $\forall \xi_1$.
Then the function $V_0(s, \xi_0) := \min_x \left\{ h(s, x, \xi_0) + \mathbb{E}\big[V_1(g(s, x, \xi_1), \xi_1) \mid \xi_0\big] \right\}$ is convex in $s$ and is monotone in the sense of $s_1 \geq s_2 \Rightarrow V_0(s_1, \xi_0) \geq V_0(s_2, \xi_0)$, $\forall \xi_0$.

PROOF 1. As $g(s, x, \xi_1)$ is convex in $s$, the following holds by definition: $g(\lambda s_1 + (1-\lambda)s_2, x, \xi_1) \leq \lambda g(s_1, x, \xi_1) + (1-\lambda) g(s_2, x, \xi_1)$, $\forall \lambda \in [0,1]$. Further, as the function $V(y, \xi_1)$ is monotonically increasing and convex, we claim that $V(g(\lambda s_1 + (1-\lambda)s_2, x, \xi_1), \xi_1) \leq V(\lambda g(s_1, x, \xi_1) + (1-\lambda) g(s_2, x, \xi_1), \xi_1) \leq \lambda V(g(s_1, x, \xi_1), \xi_1) + (1-\lambda) V(g(s_2, x, \xi_1), \xi_1)$. The same result holds with respect to the variable $x$ and, thus, the function $V_1(g(s, x, \xi_1), \xi_1)$ is jointly convex in $(s, x)$. As the sum of two convex functions, the function $h(s, x, \xi_0) + \mathbb{E}\big[V_1(g(s, x, \xi_1), \xi_1) \mid \xi_0\big]$ is also jointly convex in $(s, x)$.

In order to prove the monotonicity result, we note that $y_1 \geq y_2 \Rightarrow V_1(y_1, \xi_1) \geq V_1(y_2, \xi_1)$ and $s_1 \geq s_2 \Rightarrow g(s_1, x, \xi_1) \geq g(s_2, x, \xi_1)$, $\forall x, \xi_1$, which implies $s_1 \geq s_2 \Rightarrow V_1(g(s_1, x, \xi_1), \xi_1) \geq V_1(g(s_2, x, \xi_1), \xi_1)$, preserving the monotonicity of the function $V_1(g(s, x, \xi_1), \xi_1)$ in $s$. As the minimized sum of two monotone functions, the function $V_0(s, \xi_0)$ is also monotone in $s$ in the sense of $s_1 \geq s_2 \Rightarrow V_0(s_1, \xi_0) \geq V_0(s_2, \xi_0)$, $\forall \xi_0$. □

Analogically, one can prove a statement for monotonically decreasing functions. Overall, one can guarantee the global solution of a convex multi-stage problem with monotone functions $h_t(s_t, x_t, \xi_t)$ and $g_t(s_t, x_t, \xi_{t+1})$, by requiring stage-wise monotonicity as a constraint in the value function approximation.

### 3.1. Convex and monotone value function approximation

Lemma 1 provides the ground for the approximation of the optimal value function by a convex and monotone interpolation prior to the solution of the corresponding optimization problem. In this article, we use the QP (quadratic programming) approximation of the value function in $s_t$, i.e., $\hat{V}_t(s_t, \xi^{*i}_t) = s_t{}^T A_i s_t + 2 b_i^T s_t + c_i$, where $A_i$, $b_i$ and $c_i$ are estimated by fitting a convex and monotone function $\hat{V}_t(s_t, \xi^{*i}_t)$ to the points $\hat{V}_t(\hat{s}^k_t, \xi^{*i}_t)$,

**Algorithm 1** Dynamic programming with the linear interpolation of Fréchet quantizers.

Generate trajectories for the endogenous variable $\{\hat{s}_t^k\}_{t=1}^{T-1}$, $\forall k = 1,...,K$;
Quantize the Fréchet distribution at the stage $t = 1$ by minimizing distances (8) and (9). Interpolate the quantizers using the recusion (15);
**for** $t = T,...,1$ **do**
  **if** $t == T$ **then**
    Compute $\hat{V}_T(\hat{s}_T^k, \xi_T^{*i})$, $\forall i,k$ by solving the dynamic optimization problem (16) at stage $T$ with $\hat{V}_{T+1}(\cdot) = 0$;
  **else if** $0 \leq t \leq T-1$ **then**
    Evaluate $s_{t+1}^j = g_t(\hat{s}_t^k, x_t, \xi_{t+1}^{*j})$, $\forall j \in \mathscr{L}_{t+1}^i$;
    Interpolate $\hat{V}_{t+1}(s_{t+1}^j, \xi_{t+1}^{*j})$ by quadratic approximation (17);
    Solve the problem (16) at the stage $t$.
  **end if**
**end for**

$\forall k = 1,...,K$. The parameters can be computed efficiently by solving the following semidefinite program $\forall i = 1,...,n_t$:

$$\min_{A_i,b_i,c_i} \left\{ \sum_{k=1}^{K} \left[ (\hat{s}_t^k)^T A_i \hat{s}_t^k + 2b_i^T \hat{s}_t^k + c_i - \hat{V}_t(\hat{s}_t^k, \xi_t^{*i}) \right]^2 \right\}, \tag{17}$$
$$\text{subject to } A_i \in \mathbb{S}^{r_3},\ b_i \in \mathbb{R}^{r_3}, c_i \in \mathbb{R}$$
$$A_i \succeq 0, \quad A_i \hat{s}_t^k + b_i \geq 0,\ \forall k = 1,...,K,$$

where $\mathbb{S}^{r_3}$ is a set of symmetric matrices and $A_i \succeq 0$ denotes that the matrix $A_i$ is positive semidefinite. Note that one would not require monotonicity or convexity conditions in case of linear programming (i.e., if functions $h_t(s_t, x_t, \xi_t)$, $g_t(s_t, x_t, \xi_{t+1})$ and $V_{t+1}(s_{t+1}, \xi_{t+1})$ are linear in $s_t$ and $x_t$). Indeed, linearity conditions are a special case of requirements of Lemma 1 and they are recursively preserved in the dynamic programming. In the linear case, the interpolation of the value function can be stated in the following form:

$$\min_{b_i,c_i} \left\{ \sum_{k=1}^{K} \left[ b_i^T \hat{s}_t^k + c_i - \hat{V}_t(\hat{s}_t^k, \xi_t^{*i}) \right]^2 \right\},$$
$$\text{subject to } b_i \in \mathbb{R}^{r_3}, c_i \in \mathbb{R}.$$

Algorithm 1 describes the overall dynamic programming procedure. Note that the optimal supporting points are computed at the stage $t = 1$ only, while the estimation procedure from Section 2.2 is used at all later stages. This clearly reduces the complexity of the scenario tree quantization, which only requires the optimal transport problem to be solved at the stage $t = 1$. The time complexity of the dynamic programming, however, depends on the number of optimization problems (16) solved at each stage. For large *PNL* threshold, the number of problems is equal to $nKT$ if the sample size is large enough to assume $\hat{x}_{s,0}^i \to_{N\to\infty} \hat{x}_{s,t-1}^i \to_{N\to\infty} \hat{x}_{s,t}^i \to_{N\to\infty} \widetilde{x}_s^i$, in which case the problem becomes as follows:

$$\hat{V}_t(\hat{s}_t^k, \xi_t^{*i}) = \left\{ h_t(\hat{s}_t^k, x_t, \xi_1^{*i}) + \sum_{j \in \mathscr{L}_{t+1}^i} \left[ \hat{V}_{t+1}(s_{t+1}^j, \xi_1^{*j}) \right] p_{t+1}^{*j} \right\}, \tag{18}$$
$$\text{subject to } x_t \in \mathbb{X}_t,\ x_t \lhd \mathscr{F}_t,\ s_{t+1}^j = g_t(\hat{s}_t^k, x_t, \xi_1^{*j}),\ \forall j.$$

Also, homogeneous functions are often used to describe the final cost of natural disasters, e.g., flood events. Bayes Business School report 2023 [52] illustrates the homogeneous effect of flood risk on house prices and the consistency of the findings over the last five years. Furthermore, the study of Prahl [48] provides an economic assessment of the impacts of storm surges and sea-level rise using homogeneous functions to estimate monetary damages in coastal cities in Europe. If the functions $h_t(\cdot)$ and $V_{t+1}(\cdot)$ are homogeneous in $\xi_t^{*i}$ and $\hat{x}^i_{s,t-1} \approx \hat{x}^i_{s,t}$ (which is the case for more distant time periods $t \leq T$ in the planning horizon), the problem can be written as follows:

$$\hat{V}_t(\hat{s}_t^k, \xi_t^{*i}) = \frac{\hat{x}^i_{s,t}}{\hat{x}^i_{s,0}} \min_{x_t} \left\{ h_t(\hat{s}_t^k, x_t, \xi_1^{*i}) + \sum_{j \in \mathscr{L}^i_{t+1}} \left[ \hat{V}_{t+1}(s^j_{t+1}, \xi_1^{*j}) \right] p^{*j}_{t+1} \right\}, \tag{19}$$
$$\text{subject to } x_t \in \mathbb{X}_t,\ x_t \lhd \mathscr{F}_t,\ s^j_{t+1} = g_t(\hat{s}_t^k, x_t, \frac{\hat{x}^i_{s,t}}{\hat{x}^i_{s,0}} \xi_1^{*j}),\ \forall j.$$

In the next section, we apply the developed methods to the solution of a flood risk management problem on a governmental level in Austria.

## 4. Measuring and managing flood risk in Austria

An essential part of stresses and risks for societies and environments is imposed by natural hazards, which have low probability of occurrence but very high impact (Belval et al. [6], Timonina et al. [57], Van der Knijff [32]). Clearly, the frequency and the severity of such events is subject to climate change, which “will never stop without our intervention,” according to the host of the climate summit COP27. Therefore, the research devoted to finding optimal strategies for risk management of natural disasters is motivated by different needs of people on international, national and local policy levels. In this section, we focus on flood risk in Austria, considering it from multi-period and multi-regional points of view. We analyze economic risks imposed by flood events and develop a strategy which promotes the adaptation, resilience and resistance of societies to catastrophes and contributes to a decrease of risk and vulnerability. Due to the uncertainty about place and time of occurrence of flood events and in order to find the best strategies for mitigation and reconstruction after them, we assess the interregional uncertainty using the structured coupling approach as of Timonina et al. [57]. Constructing the optimal strategy, we especially take cooperation measures such as governmental insurance mechanisms and infrastructure investments into account. As the UN Secretary-General at the COP27 Climate Implementation Summit put it: ”Humanity has a choice: cooperate or perish.” Formulating a problem for risk management of catastrophic events in terms of a multi-stage stochastic optimization problem with disaster losses described by Fréchet distributions, we find optimal pre-event and post-event strategies using the developed solution method on scenario trees (Algorithm 1).

### 4.1. Multi-regional risk analysis

We consider flood events in Austria as an example of rare but damaging events. Figure 4 shows European and, particularly, Austrian river basins subject to flood risk. In Figure 4a, one can observe the dense network structure of river basins in Europe, while Figure 4b shows separate river basins in Austria.

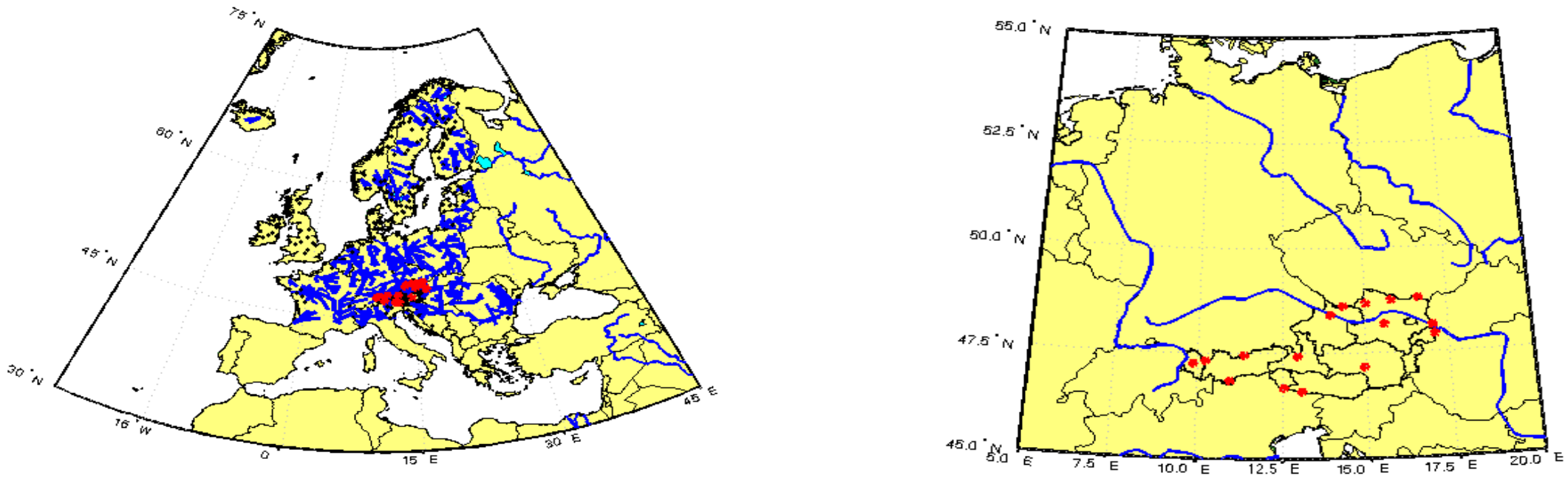

a. Europe. b. Austria.

**Figure 4 River basins in Europe and in Austria subject to flood risk.**

To estimate the joint flood risk in any combination of these basins, it is not enough to sum up their losses and to estimate the distribution function. This is due to potential non-linear interdependencies between regions, which can be more interdependent in case of severe events and less interdependent in case of more frequent events. To estimate the joint risk of flood events in a region (e.g., on local, national or regional levels), one needs to use vine copulae as of Bedford [4] and Kurowicka [33] or their approximations. For vine copulae, one would be required to estimate $(B-1)!$ copula functions, where $B$ is the number of interdependent basins. For efficiency in the presence of a high number of basins, we use a quick structured coupling approach proposed by Timonina et al. [57]. The benefit of this method is that it allows to avoid underestimation of risks without the use of robust optimization. The marginal losses in each basin are simulated using the LISFLOOD hydrological model and an economic damage model as in the works of Van der Knijff, Younis and De Roo [32], as well as Rojas et al. [49].

Coupling marginal losses in such a way that the large-scale probability distribution is not underestimated and fits the multi-regional data on losses, we find a non-linear dependency function (i.e., copula) that hierarchically transforms pairs of marginal probability distributions into the joint probability distribution taking all the regional interdependencies into account (Timonina et al. [57]). The order of coupling is directly based on the available data on the river structure. Overall, we estimate national-scale probability loss distributions for Austria using Flipped Clayton, Gumbel and Frank copulae (see Table 3 for 5-, 10-,..., 1000-year events and Probabilities of No Loss, i.e., PNL).

| Year events | | 5 | 10 | 20 | 50 | 100 | 250 | 500 | 1000 |
|---|---|---|---|---|---|---|---|---|---|
| Cum. probabilities | | 0.8 | 0.9 | 0.95 | 0.98 | 0.99 | 0.996 | 0.998 | 0.999 |
| Fl. Clayton ($L$) | PNL=0.68 | 0.000 | 0.030 | 0.855 | 2.926 | 8.017 | 13.810 | 16.756 | 17.761 |
| Gumbel ($L$) | PNL=0.67 | 0.000 | 0.067 | 0.847 | 2.953 | 7.975 | 13.107 | 16.340 | 17.736 |
| Frank ($L$) | PNL=0.63 | 0.000 | 0.395 | 1.147 | 2.894 | 7.697 | 10.213 | 10.689 | 10.689 |

**Table 1 Total losses in Austria for 2030 in EUR bln.**

Note that Flipped Clayton and Gumbel copulae assume stronger dependencies for higher rather than smaller losses. Differently, the Frank copula considers the dependency between lower and upper tails as symmetrical. Thus, due to their asymmetric structure, the Flipped Clayton copula or the Gumbel copula are preferred over the Frank copula for modeling interdependencies in flood losses. The Frank copula would give the same weights to events of higher and lower frequencies (see also Appendix A.2). Figure 5a demonstrates the national-scale loss distribution estimates obtained via the use of Flipped Clayton Copula for two past and two future years, i.e., 1995, 2015, 2030 and 2050. The continuous fits of Fréchet, Weibull and Gumbel distributions, as well as the average value-at-risk and the expected values for year 2030 are demonstrated in Figure 5b. Note that the climate change influence is pronounced in years 1995 and 2015, implying an increase in the frequency of extreme events. By this, the distributions for years 2030 and 2050, whose estimates do not account for new major events, can be over-optimistic. In general, when a new disaster event occurs, an update of the loss distribution parameters is required as described in Section 2.2. As one can see in Figure 5a, the Weibull and Gumbel distributions underestimate the risk of low-probability events in Austria. Furthermore, the Gumbel distribution highly overestimates the risk of more frequent events.

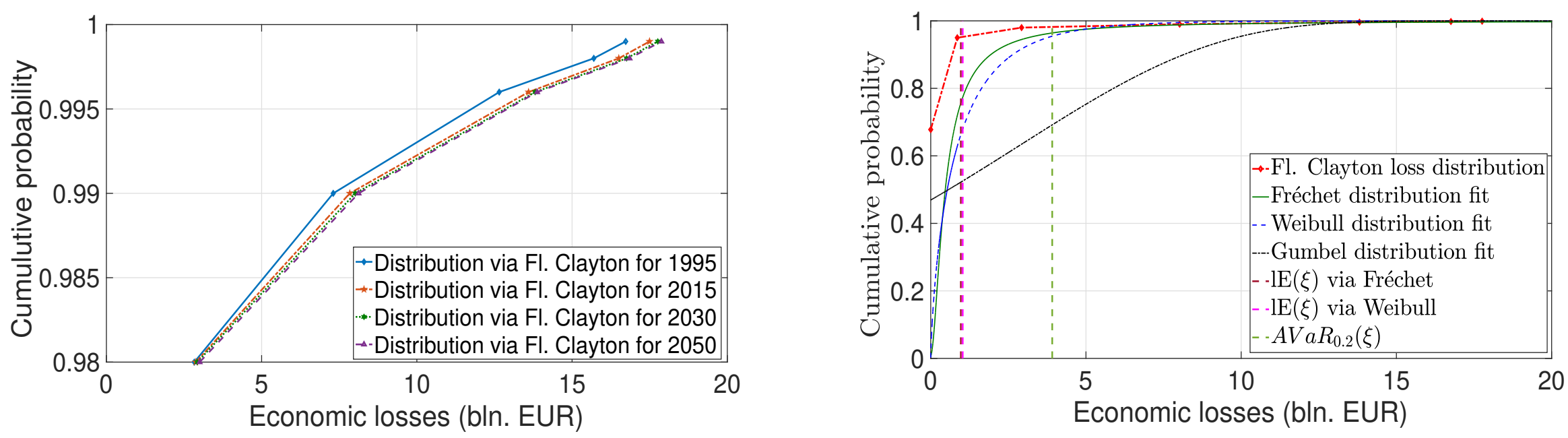

a. Loss distributions between 1995 and 2050.

b. Heavy-tailed fits for year 2030.

**Figure 5 Estimated flood risk in Austria.**

We use the estimated Fréchet distribution because it represents a continuous heavy-tailed distribution function, which avoids underestimation of losses after low-probability events. Moreover, it gives the lowest Kantorovich-Wasserstein metric between the coupled distribution and the fitted one (see Kantorovich [30]). To represent the flood risk in a *multi-period* environment, we approximate the stochastic process $\xi = (\xi_1, \xi_2, ..., \xi_T)$, which describes random losses $\xi_t$ after flood events in each period $\forall t = 1, ..., T$, via a finitely-valued scenario tree. Here, Figure 6a demonstrates the fast quantization method from Section 2.2; Figure 6b shows the method with varying $\lambda$ as of Gumbel [24]; Figure 6c stands for the combination of the two. Overall, the risk of flood events is represented by a scenario tree based on the optimally-quantized Fréchet distribution. Note that high risks are not represented in Figure 6a. This can be suitable for countries with stable economies (e.g., Austria) and happens due to the fact that the shape parameter is considered constant throughout the planning horizon. Oppositely, Figures 6b and c demonstrate patterns with more pronounced heavy-tailed quantizers. Note that the scenario tree in Figure 6c is more efficient to generate due to the use of a fast approach for $100\% \cdot \mathrm{PNL}$ of nodes. Accounting for low-probability high-impact floods, we use the multi-period risk representation as in Figure 6c to formulate and solve a multi-stage optimization problem for risk management of flood events.

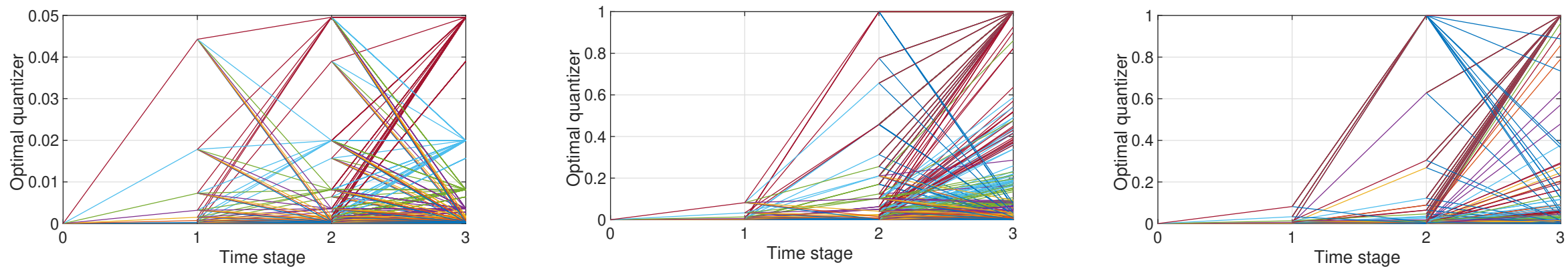

a. Tree with constant $\lambda$. b. Tree with varying $\lambda$. c. Combined method.

**Figure 6 Quantization of the Fréchet distribution.**

## 4.2. Multi-period risk management

We consider a government which may lose a part of its capital $S_t$ at any time $t = 1,\ldots,T$ because of random natural hazard events with uncertain relative economic loss $\xi_t$, i.e., $\xi_t \in [0,1]$. As a result of this loss, the country would face a drop in GDP at the end of the year. We focus on two major components of the *government expenditures*, which are i) *gross capital formation*, that can be classified as *government investment*, and ii) *government consumption*, whose purpose is to satisfy the individual or collective needs of the community. These two types of government expenditures together play a crucial role for the *gross domestic product*. Importantly, the amounts which should be allocated to the government investment and to the government consumption are subject to the decision of the government in a setup without natural disasters. However, in a setup in which a part of the capital might be destroyed due to natural disasters, the problem becomes more complex and leads to the necessity to allocate a part of the budget to risk mitigation strategies (see Figure 7).

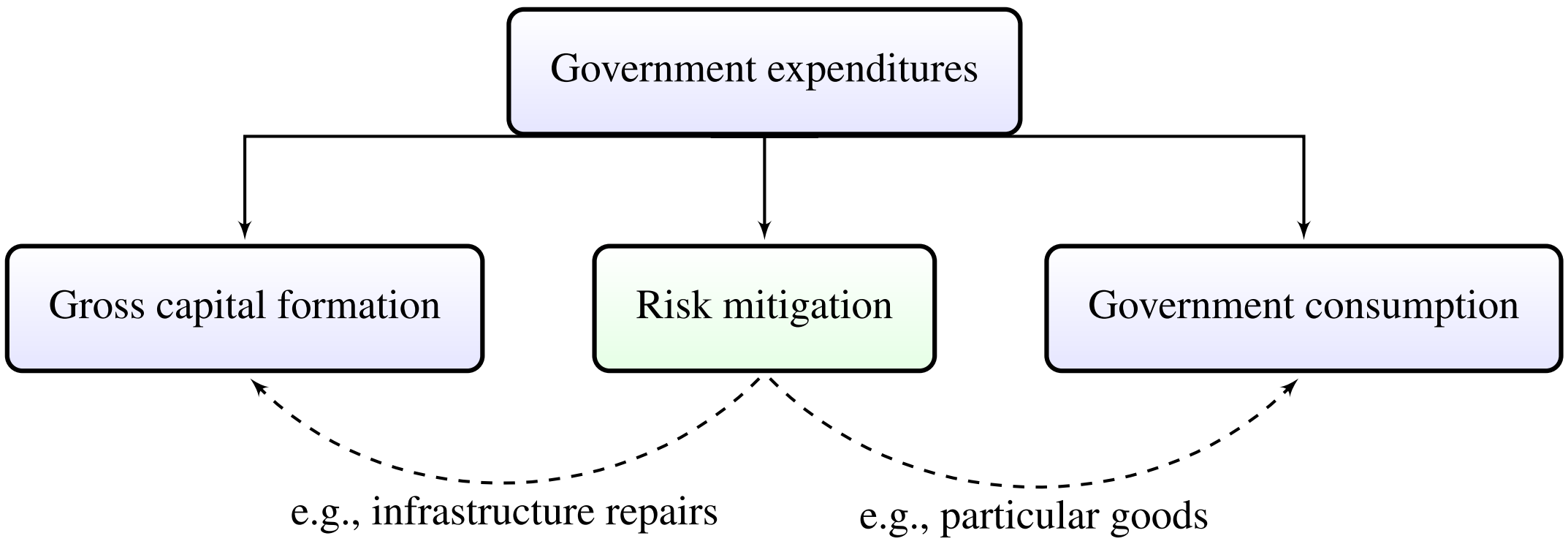

**Figure 7 Structure of government expenditures**

Therefore, we consider a stochastic optimization problem to decide on the optimal allocation of the available budget $B_t$ to i) the investment $x_t$, which influences the capital formation, ii) the government consumption $c_t$ and iii) the risk mitigation strategies $z_t$. Incorporating risk mitigation strategies into the analysis, we assume that a government insurance scheme against natural disasters is available in the country and that the decision about the amount $z_t$, $\forall t = 0\ldots,T-1$ needs to be made. The insurance premium depends on the expected value of the relative loss $\xi_t$, $\forall t = 1,\ldots,T$ and is equal to $\pi(\mathbb{E}(\xi_t)) = (1+V)\mathbb{E}(\xi_t)$, where $V$ is a constant insurance load. All decisions are taken one stage prior to the realization of the loss $\xi_t$. As we consider the flood mitigation problem in a multi-period setup, we generalize the framework proposed by van Danzing [15] to cases where the flood loss distribution does not stay constant over time. We also do not optimize over the height of dikes (see Eijgenraam et al. [18], van Danzing [15] for the works on flood prevention in the Netherlands), but focus on economically efficient flood mitigation strategies, allowing to reduce losses in case of events (see Eijgenraam et al. [19]).

Dependent on the goals of the government, different decision-making problems can be stated in terms of objective functions and the constraints. We consider a multi-stage model, which describes the decision-making problem in terms of relative capital loss $\xi_t$, $\forall t = 1,...,T$, while the GDP is being modeled in line with the classical Cobb-Douglas production function with constant productivity and maximal weight of capital rather than labor. The available budget is a constant part of the GDP in this setting. Hence, $B_t = \alpha S_t$, where $S_t$ is the governmental capital at stages $t = 0,...,T$ and $\alpha$ is a constant term from the interval $[0,\ 1]$.

In terms of governmental objective, we consider the multi-stage stochastic optimization program, where $\beta$, $\delta$, $\rho$ are constant factors from the interval $[0,\ 1]$ and $S_0 = 322.56$ bln. EUR:

$$\max_{x_t,c_t,z_t} \mathbb{E}\left[(1-\beta)\sum_{t=0}^{T-1}\rho^{-t}u(c_t) + \beta\rho^{-T}u(S_T)\right] \tag{20}$$

subject to $x_t,\ z_t,\ c_t \geq 0,\ t = 0,...,T-1;\ S_0$ is given.

$$x_t \lhd \mathscr{F}_t,\ c_t \lhd \mathscr{F}_t,\ z_t \lhd \mathscr{F}_t,\ t = 0,...,T-1,$$

$$S_{t+1} = [(1-\delta)S_t + x_t](1-\xi_{t+1}) + z_t\xi_{t+1},\ t = 0,...,T-1,$$

$$B_t = \alpha S_t \geq x_t + c_t + \pi(\mathbb{E}(\xi_{t+1}))z_t,\ t = 0,...,T-1.$$

The complexity of the problem arises as the stochastic process $\xi = (\xi_1,...,\xi_T)$ is interdependent in time, meaning that all future random variables depend on past realizations of the process.

In the problem (20), $u(\cdot)$ is a governmental utility function of the form $u(c) = \frac{c^{1-\gamma}}{1-\gamma}$. Such utility can differentiate between risk-neutral and risk-averse risk bearing abilities dependent on the chosen parameter $\gamma \in [0,1]$ (see Hochrainer and Pflug [29] for more details on governmental risk aversion). The objective of the decision-maker is to maximize the expectation about the weighted government consumption, whose aim is to represent overall individual and collective satisfaction of the community at each period $t$, and the government capital $S_T$ at the final stage, whose purpose is to provide enough resources for the future. Note that the discounting factor $\rho$ assigns non-decreasing weights to future capital and consumption.

Solving the multi-stage problem numerically, we obtain approximations of stage-wise optimal decisions and the optimal value. For the endogenous variable $S_t$, we sample $K$ random trajectories $\{S_t^k\}_{t=1}^T$ from the interval $(0,(1-\delta+\alpha)^{t-1}S_0]$, $\forall t = 1,...,T$, where the upper limit is the capital in case of no event and the maximal investment amount $x_t = \alpha S_t$. In line with Algorithm 1, we rewrite problem (20) in the dynamic programming form, which is solved iteratively starting with the stage $T-1$ and going backwards to the stage $t = 0$. As described before, we use optimal weighting of value functions in our solution method. Importantly, this allows to avoid underprevention of flood risk in Austria. Indeed, according to Baillon et al. [3], underprevention is very often being caused by misperceived probabilities. We compare the quality of the numerical solution for constant and varying patterns in the shape parameter $\lambda$ of the Fréchet distribution. We also use a combined approach, where the probability of no loss serves as an indicator for the incorporated parameter $\lambda$ (i.e., it helps to decide if one should use a constant or a varying pattern, see Figure 8). As expected, we observe a convergence of the optimal value dependent on tree branchiness. While using the constant parameter $\lambda$ for all tree nodes may lead to the underestimation of the optimal value (Figure 8a), the combined approach demonstrates the convergence with less variability of estimates in comparison with the fully varying pattern. Clearly, the combined approach also gains in efficiency due to the use of a constant shape parameter with high probability.

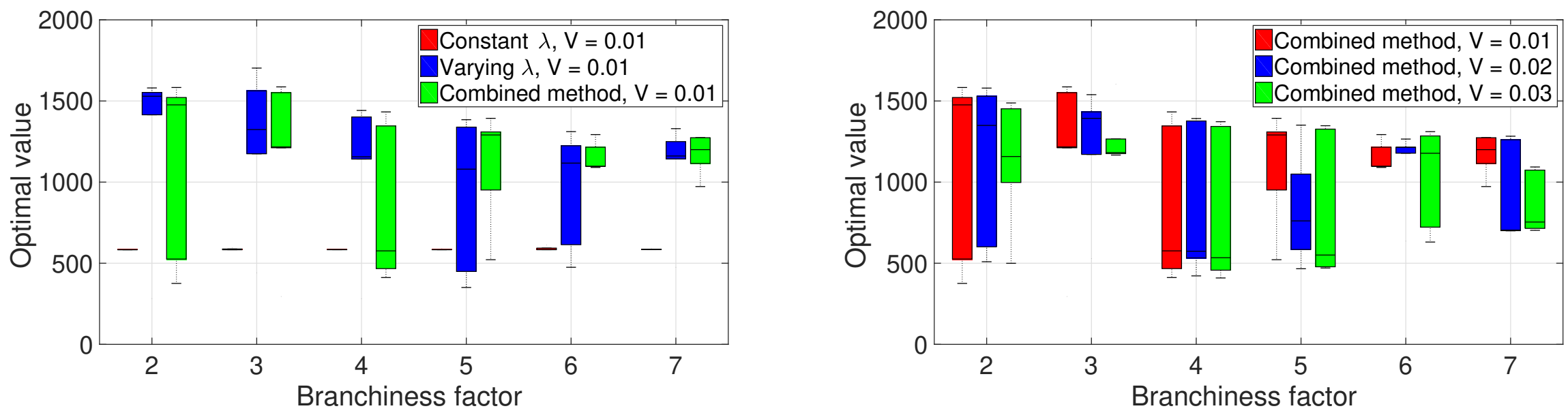

a. Comparison of methods. b. Combined approach.

**Figure 8 Optimal values for the problem (20) obtained using different patterns of the shape parameter $\lambda$.**

In Figure 9, one can observe optimal patterns of insurance and investment decisions dependent on time and the insurance load parameter $V$. Given the planning horizon of $T = 10$ years in Figure 9a, one can observe that the average insurance amount tends to decrease towards the middle of the planning horizon and to increase afterwards. This is a property of the time-varying flood risk in a region and possible underestimation of this risk by the expectation $\mathbb{E}(\xi_t)$: indeed, for many high-risk scenarios, the cost of the insurance at later stages appears to be too low. Similar behavior can be observed for the optimal investment in Figure 9b. Note that the average insurance amount is much higher than the average investment given low enough insurance costs (i.e., $0.01 \leq V \leq 0.03$). Following this, one can conclude that Austria could benefit from lost-cost insurance schemes facing the risk of low-probability high-impact floods.

Importantly, the amount of the optimal insurance at later stages increases on average until some level of the load parameter $V$ is reached (Figure 10a). This is also observed in Figure 9a for three values of the parameter $V$ and is a consequence of the insurance payoff function. Indeed, if the cost of the insurance increases, one could decrease the amount of the insurance purchased to stay on the same budget level (see the budget constraint in problem (20)). However, the payoff in case of a disaster event would also decrease in this case, followed by the drop in governmental stock and in the objective function. Thus, one would

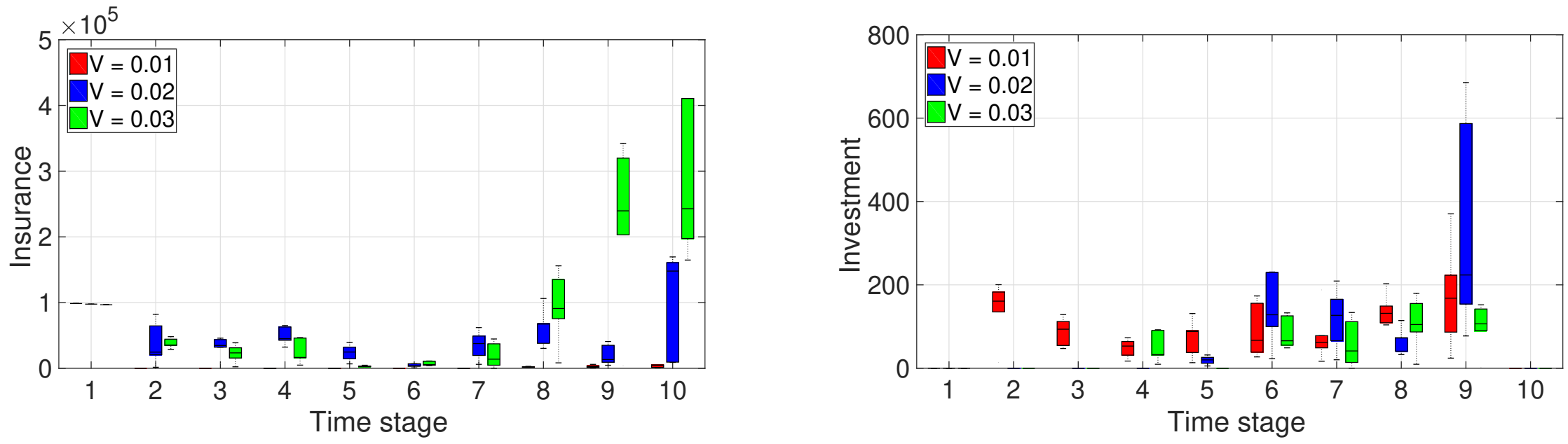

a. Insurance decision. b. Investment decision.

**Figure 9 Optimal insurance and investment decisions for different values of the insurance load $V$.**

need to increase an investment, a consumption or, again, an insurance so that the objective value does not decrease. Note that the increase in consumption levels may be suboptimal due to the non-linearity of the objective function. In this case, the insurance amount increases despite the rising costs.

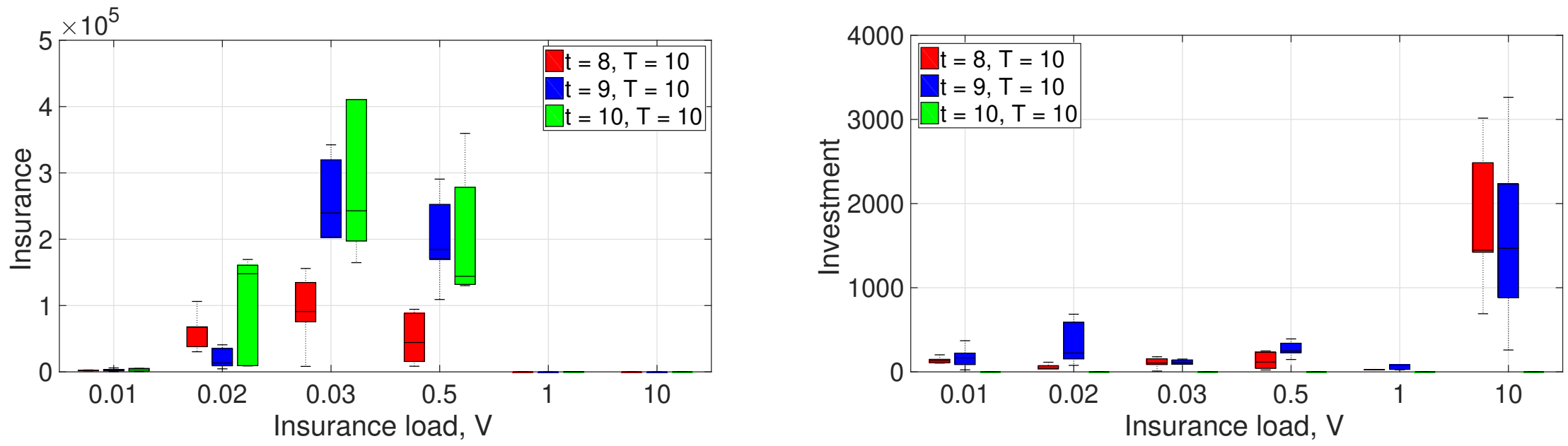

a. Lower values of parameter $V$. b. Dependency on parameter $V$.

**Figure 10 Optimal investment decision for different values of the insurance load $V$.**

Next, as we know optimal probabilities of scenarios by solving the Kantorovich-Wasserstein minimization problem, we can compute capital distribution functions given different values of the insurance load $V$. Note that this would not be possible in case of Monte-Carlo sampling as the scenarios would be equally probable. The capital distribution accounts for optimal decisions in the problem (20), i.e., this is the distribution given the maximal objective value. Clearly, the distribution would change under different (i.e., suboptimal) allocation strategies (Figure 11).

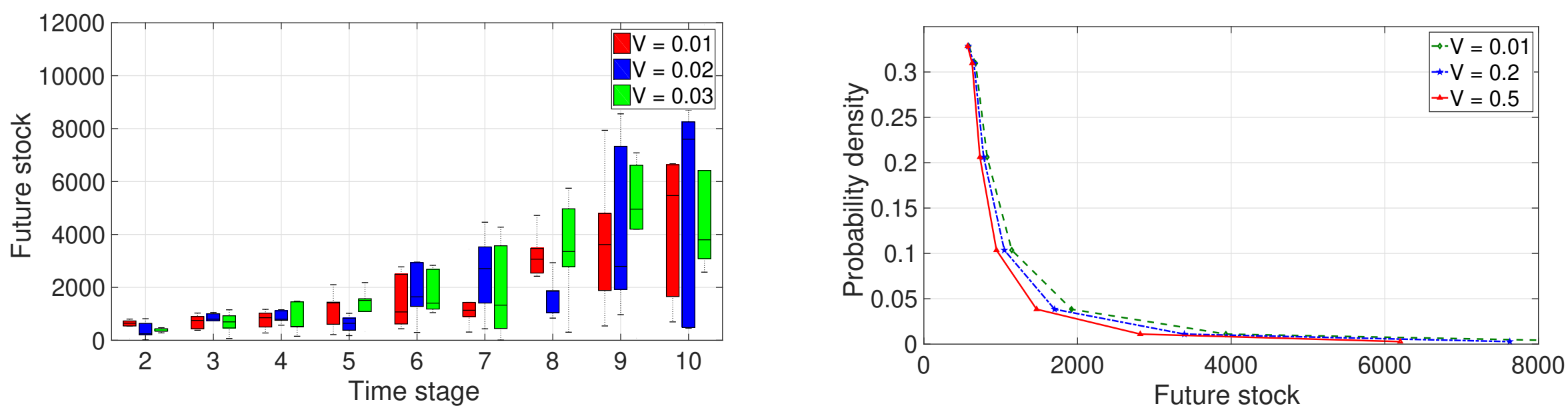

a. Stock building pattern. b. Stock distribution.

**Figure 11 Future stock given the optimal decision-making.**

In Figure 11a, one can see the development of governmental stock in time: optimal allocation guarantees the growth of Austrian stock on average. Next, in Figure 11b, one clearly observes that more expensive insurance policies (i.e., $V = 0.2$ or $V = 0.5$) lead to lower capital in $T = 10$ years in probability. These are also the policies which make the insurance purchase suboptimal (Figure 10a).

### 4.3. Robust model for risk-management of extremes

Next, suppose that the conditional probabilities used in the optimization problem (16) may deviate from their nominal values $p_t^{*i}$, $\forall i = 1, \dots, n$ up to a certain degree (see Hanasusanto et al. [25] for more details), i.e., let us consider the following ambiguity set, denoted by $\Delta(p_t^*)$:

$$\Delta(p_t^*) = \{q \in \Delta : \sum_{i=1}^{n} \frac{(p_t^{*i} - q^i)^2}{q^i} \le \theta\}, \tag{21}$$

where $p_t^* = (p_t^{*1}, \dots, p_t^{*n})$ and $\theta$ is a risk budget. Note that the set (21) simplifies to the standard simplex in a degenerate case for large enough $\theta$ (see Figure 12).

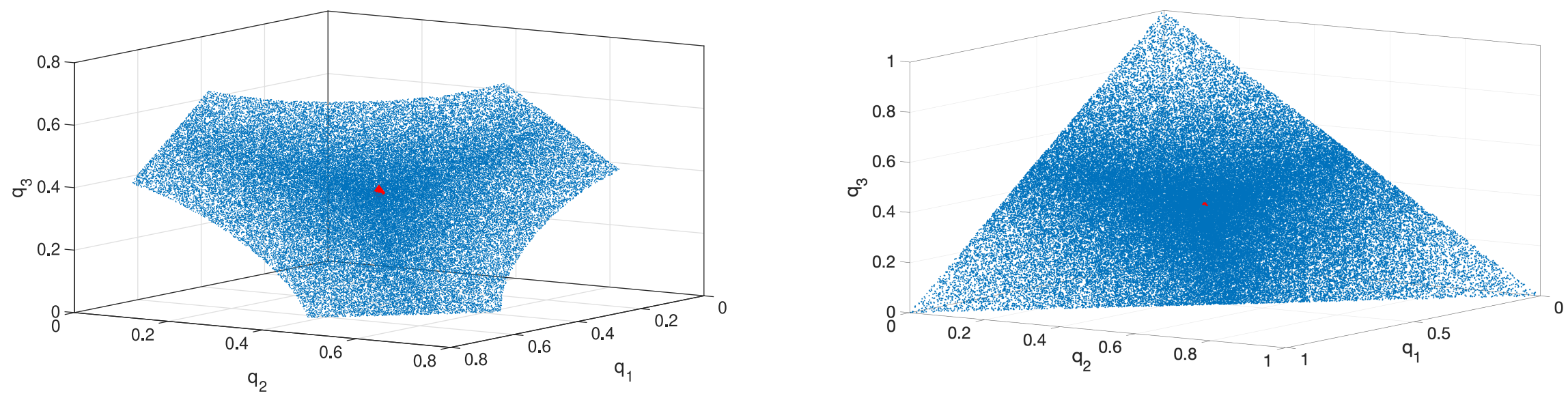

a. $\theta = 0.1$ b. $\theta = 10$

**Figure 12 Three-imensional ambiguity set for $\theta = 0.1$ and $\theta = 10$.**

With a small abuse of notations, the dynamic representation of the problem (20) can be written in the following distributionally-robust form:

$$\begin{aligned}\hat{V}_t(\hat{s}_t^k, \xi_t^{*i}) = \; & \max_{x_t} \Big\{ h_t(\hat{s}_t^k, x_t, \xi_t^{*i}) + \min_{q \in \Delta(p_t^*)} \sum_{j \in \mathcal{L}_{t+1}^i} \big[\hat{V}_{t+1}(s_{t+1}^j, \xi_{t+1}^{*j})\big] q^j \Big\}, \\ & \text{subject to } x_t \in \mathbb{X}_t,\ x_t \lhd \mathscr{F}_t,\ s_{t+1}^j = g_t(\hat{s}_t^k, x_t, \xi_{t+1}^{*j}),\ \forall j \in \mathcal{L}_{t+1}^i, \\ & \xi_{t+1}^{*j} = \frac{\hat{x}_{s,t}^i(\hat{s}_t^k)}{\hat{x}_{s,0}^i} \xi_1^{*j},\ \forall j \in \mathcal{L}_{t+1}^i,\end{aligned} \tag{22}$$

which leads to the stage-wise reformulation of the flood risk-management problem (20).

THEOREM 1. Stage-wise reformulation of the flood risk-management problem (20) in the uncertainty set (21) can be written as follows $\forall t < T$:

$$\begin{aligned}&\max_{\substack{x_t, c_t, z_t \\ s_{t+1}^i, \forall i \\ \mu_1, \mu_2, y, v}} \big[(1-\beta)\rho^{-t} u(c_t) + \mu_1 \theta - \mu_2 + 2(p_t^*, y) - 2\mu_1(p_t, \mathbb{1})\big] \\ &\text{subject to } x_t \ge 0,\ c_t \ge 0,\ z_t \ge 0,\end{aligned} \tag{23}$$

$$\mu_1 \in \mathbb{R}_{\{0,+\}},\ \mu_2 \in \mathbb{R},\ y, v \in \mathbb{R}^n,$$
$$\alpha \hat{s}_t^k \geq x_t + c_t + \pi(\xi_{t+1}) z_t,$$
$$s_{t+1}^i = [(1-\delta)\hat{s}_t^k + x_t](1-\xi_{t+1}^{*i}) + \xi_{t+1}^{*i} z_t,\ \forall i$$
$$v_i + \mu + \lambda > 0,\ v_i \leq \hat{V}_{t+1}(\hat{s}_{t+1}^k, \xi_{t+1}^{*i}),\ \sqrt{4y_i^2 + (v_i+\mu)^2} \leq 2\lambda + v_i + \mu,\ \forall i.$$

PROOF 2. Consider the following problem:

$$\min_{q \in \mathbb{R}_+^N} v^T q, \text{ subject to } \mathbb{1}^T q = 1,\ \sum_{i=1}^N \frac{(p_i - q_i)^2}{q_i} \leq \theta,$$

where $\mathbb{1}$ is the all-ones vector and $v \in \mathbb{R}^n$ is a vector with entries $\hat{V}_{t+1}(\hat{s}_{t+1}^k, \xi_{t+1}^{*i}),\ i = 1, ..., n$.
The Lagrangian of the problem is defined as

$$\mathscr{L}(q, \mu_1, \mu_2) = v^T q + \mu_1 \left( \sum_{i=1}^n \frac{(p_t^{*i} - q_i)^2}{q_i} - \theta \right) + \mu_2 \left( \mathbb{1}^T q - 1 \right),$$

where $\mu_1 \geq 0$ and $\mu_2 \in \mathbb{R}$ are the dual multipliers. By strong duality, the following holds $\min_{q \in \mathbb{R}_+^n} v^T q = \min_{q \in \mathbb{R}_+^n} \max_{\substack{\mu_1 \in \mathbb{R}_+ \\ \mu_2 \in \mathbb{R}}} \mathscr{L}(q, \mu_1, \mu_2) = \max_{\substack{\mu_1 \in \mathbb{R}_+ \\ \mu_2 \in \mathbb{R}}} \min_{q \in \mathbb{R}_+^n} \mathscr{L}(q, \mu_1, \mu_2)$. Considering the second equality and noting that $\mathscr{L}(q, \mu_1, \mu_2) = -\mu_1 \theta - \mu_2 + \sum_{i=1}^n \left( (v_i + \mu_2) q_i + \mu_1 \frac{(p_t^{*i} - q_i)^2}{q_i} \right)$, we can solve the inner minimization problem in a closed form, i.e., obtaining the optimal solution for $\mu_1 + \mu_2 + v_i > 0$, which yields $q_i = p_t^{*i} \sqrt{\frac{\mu_1}{\mu_1 + \mu_2 + v_i}},\ \forall i = 1, ..., n$ at optimality. Note, that this solution is the global minimum as the objective function is convex in $q_i$. Now, we can state the following equivalence:

$$\min_{q \in \mathbb{R}_+^n} v^T q = \max_{\substack{\mu_1 \in \mathbb{R}_+ \\ \mu_2 \in \mathbb{R}}} \left[ -\mu_1 \theta - \mu_2 + \sum_{i=1}^n \left( 2 p_t^{*i} \sqrt{\mu_1 (v_i + \mu_1 + \mu_2)} - 2 p_t^{*i} \mu_1 \right) \right]$$

if $\mu_1 + \mu_2 + v_i > 0,\ \forall i = 1, ..., n$. Denoting $\sqrt{\mu_1 (v_i + \mu_1 + \mu_2)}$ by $y_i$, where $y = (y_1, ..., y_n)$, we can conclude that the problem (22) has the following convex optimization form:

$$\begin{aligned} \hat{V}_t(\hat{s}_t^k, \xi_t^{*i}) = \ & \max_{x_t} \left\{ h_t(\hat{s}_t^k, x_t, \xi_t^{*i}) - \mu_1 \theta - \mu_2 + 2(p_t^*, y) - 2\mu_1 (p_t^*, \mathbb{1}) \right\}, \\ & \text{subject to } x_t \in \mathbb{X}_t,\ x_t \triangleleft \mathscr{F}_t,\ s_{t+1}^j = g_t(\hat{s}_t^k, x_t, \xi_{t+1}^{*j}),\ \forall j \in \mathscr{L}_{t+1}^i, \\ & \xi_{t+1}^{*j} = \frac{\hat{x}_{s,t}^i(\hat{s}_t^k)}{\hat{x}_{s,0}^i} \xi_1^{*j},\ \forall j \in \mathscr{L}_{t+1}^i, \\ & v_i + \mu_1 + \mu_2 > 0,\ v_i \leq \hat{V}_{t+1}(\hat{s}_{t+1}^k, \xi_{t+1}^{*i}),\ \forall i \\ & \sqrt{4y_i^2 + (v_i + \mu_2)^2} \leq 2\mu_1 + v_i + \mu_2,\ \forall i, \end{aligned} \tag{24}$$

from which the statement of the Theorem follows directly. □

We use the robust reformulation (23) in order to test the sensitivity of the optimal solution to an increase in the risk budget $\theta$, i.e., to the increase in the uncertainty about the distribution function. As expected, the optimal value decreases as the risk budget increases (Figure 13b), while the investment and the consumption decisions deviate in opposite from each other directions (Figure 13b). Importantly, the sensitivity is stronger for higher values of the risk budget $\theta$. As our flood loss distribution can be seen as a worst-case distribution for extremes by construction in Section 2.2, the use of lower values of $\theta$ is motivated.

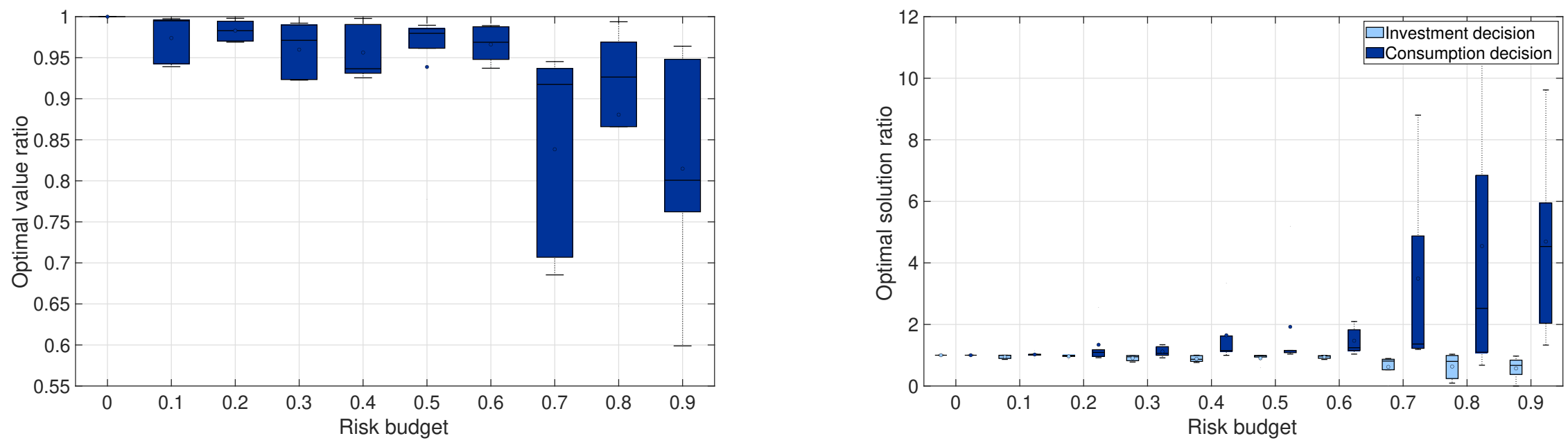

a. Optimal value. b. Optimal solution.

**Figure 13 Sensitivity of the optimal solution to an increase in the risk budget $\theta$.**

## 5. Conclusion

The focus of this article is on numerical methods for the solution of multi-stage stochastic optimization problems with uncertainty described by the Fréchet distribution. The distribution may change in time conditionally on past realizations. We propose a highly efficient *scenario grouping procedure* based on changing medians, allowing for the classification of realizations according to their impact.

Our primary focus is on applications in flood risk management, where we integrate dynamic programming techniques with a stage-wise optimal quantization method for flood loss distributions. With the current worsening climate crisis, research that focuses on optimal strategies for risk management in the face of increasing climate change-induced natural disasters will become ever more relevant. Although the methods devised in our research are aimed at flood risk management, they are also applicable across a wide range of fields characterized by Fréchet-type uncertainty. These encompass areas like flood and wind risk management, optimization of interfacial damage in microelectronic packages, and analysis of material properties of particles in engineering applications (Abbas et al. [1], Harlow et al. [27]).

Combining dynamic programming with the stage-wise optimal quantization of Fréchet distributions, we guarantee the global optimum under convexity and monotonicity conditions on the value function. In order to guarantee that these conditions hold, we introduce a monotonicity constraint in the convex optimization problem required for the value function approximation. Nevertheless, this constraint can be omitted if the value function is monotone by definition, e.g., linear or exponential.

Possible generalizations of our algorithms include data-driven methods avoiding estimation of continuous distribution functions. Such methods are usually very efficient but do not provide guarantees on the global optimum. Thus, data-driven schemes with guarantees for a global solution and quantization algorithms with utility-based metrics (e.g., Abdellaoui et al. [2], Petracou et al. [39]) are of high interest for future research and lead to a variety of applications (Bertsimas et al. [10], Piccialli et al. [45], Deng et al. [16]). Also, as the most favorable outcomes in this article are achieved when the objective function is homogeneous in the random variable, possible generalizations provide a clear avenue for future research.

# 6. Acknowledgements

We thank the Enterprise for Society Center (an UNIL, IMD and EPFL collaboration) for funding our research under the grant number E4S-3-3-03.

# Electronic Companion to "A Quick Estimation of Fréchet Quantizers for a Dynamic Solution to Flood Risk Management Problems"

## A.1. Scenario trees and the approximation quality

DEFINITION A1. A stochastic process $\xi = (\xi_1, ..., \xi_T)$ is called *a tree process* (see Pflug and Pichler [42, 44], Römisch [50]), if $\sigma(\xi_1), \sigma(\xi_2), ..., \sigma(\xi_T)$ is a filtration[1].

The *history process* $(\xi^1, \xi^2, ..., \xi^T)$ of the stochastic process $\xi$ is a tree process by definition, as soon as $\xi^1 = \xi_1$, $\xi^2 = (\xi_1, \xi_2)$,..., $\xi^T = (\xi_1, \xi_2, ..., \xi_T)$. If the stochastic process is finitely valued (e.g., $\widetilde{\xi} = (\widetilde{\xi}_1, ..., \widetilde{\xi}_T)$), it is representable as *a finitely valued scenario tree*[2] (Timonina [56]). In order to work with general tree structures, we let $n_t$, $\forall t = 1, ..., T$ be the total number of scenarios at the stage $t$ and $n_t^i$ ($\forall i = 1, ..., n_{t-1}, \forall t = 2, ..., T$) be the number of quantizers corresponding to the $n_{t-1}$ conditional distributions sitting at the stage $t$ ($\forall t = 2, ..., T$). We note that $n_t = \sum_{i=1}^{n_{t-1}} n_t^i, \forall t > 1$.

DEFINITION A2. Consider a finitely valued stochastic process $\widetilde{\xi} = (\widetilde{\xi}_1, ..., \widetilde{\xi}_T)$ represented by the tree with the same number of successors $b_t$ at each node at a given stage $t$, i.e., $n_t^i = b_t$, $\forall i$. The vector $b = (b_1, ..., b_T)$ is *a branchiness vector* of the tree (e.g., Timonina [56]). Values $b_1, b_2, ..., b_T$ are the *branchiness factors* at corresponding stages.

Figure A2 shows two randomly sampled trees with different branchiness factors. The tree on the left-hand side is a *binary tree* and, therefore, its branchiness vector is $b = [2\ 2\ 2]$; the tree on the right-hand side is a *ternary tree* and, hence, its branchiness vector is $b = [3\ 3\ 3]$. Note that both scenario trees are not optimal as they do not guarantee the minimization of Kantorovich-Wasserstein distances at stages $t = 1, 2, 3$.

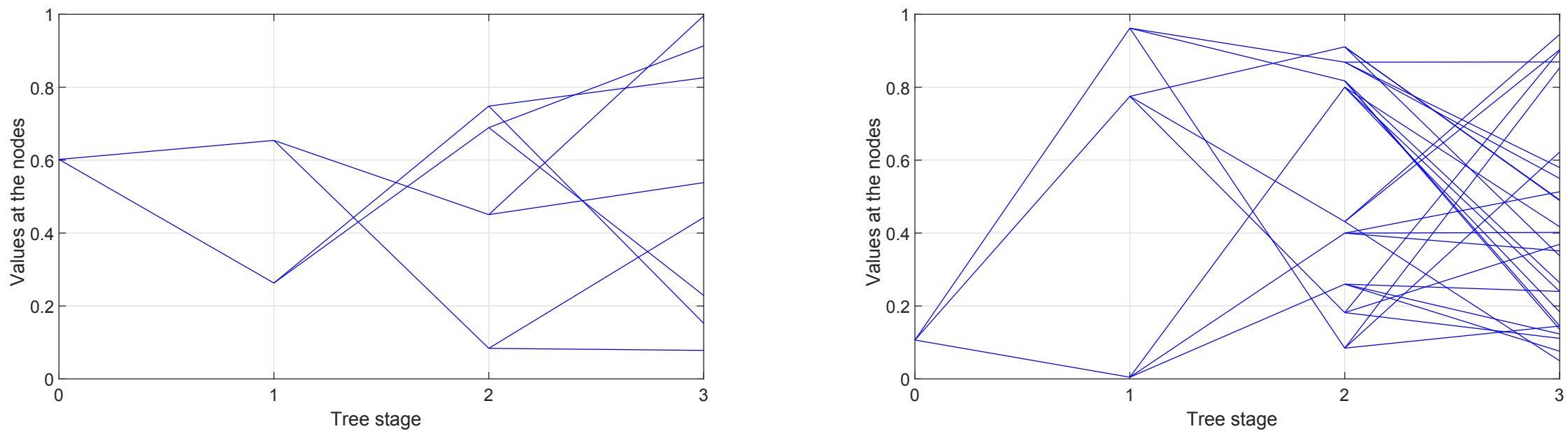

**Figure A2** **Scenario trees with different branchiness:** $b = [2, 2, 2]$ **and** $b = [3, 3, 3]$ **correspondingly.**

Both trees demonstrate the univariate case (i.e., $\widetilde{\xi}_t$ is one-dimensional $\forall t$) and, therefore, values sitting on the nodes are shown in Figure A2. In case of multidimensionality of the stochastic process $\widetilde{\xi}$, multidimensional vectors would correspond to each node of the tree and, hence, graphical representation as in Figure A2 would not be possible. Clearly, scenario trees in Figure A2 are not optimized.

## A.2. Total losses in Austria via different models

[1] Given a measurable space $(\Omega, \mathscr{F})$, a filtration is an increasing sequence of $\sigma$-algebras $\{\mathscr{F}_t\}$, $t \geq 0$ with $\mathscr{F}_t \subseteq \mathscr{F}$ such that: $t_1 \leq t_2 \Longrightarrow \mathscr{F}_{t_1} \subseteq \mathscr{F}_{t_2}$. In our case, $\sigma(\nu)$ is the $\sigma$-algebra generated by the random variable $\nu$.

[2] The scenario tree, which represents the finitely valued stochastic process $(\widetilde{\xi}_1, ..., \widetilde{\xi}_T)$, is called *a finitely valued scenario tree*.

| Year events | | 5 | 10 | 20 | 50 | 100 | 250 | 500 | 1000 |
|---|---|---|---|---|---|---|---|---|---|
| Cum. probabilities | | 0.8 | 0.9 | 0.95 | 0.98 | 0.99 | 0.996 | 0.998 | 0.999 |
| Fl. Clayton ($L$) | PNL=0.68 | 0.000 | 0.030 | 0.843 | 2.878 | 7.839 | 13.599 | 16.505 | 17.495 |
| Gumbel ($L$) | PNL=0.67 | 0.000 | 0.067 | 0.835 | 2.903 | 7.796 | 12.897 | 16.072 | 17.468 |
| Frank ($L$) | PNL=0.62 | 0.000 | 0.390 | 1.140 | 2.844 | 7.528 | 10.053 | 10.523 | 10.523 |

**Table 2** **Total losses in Austria for 2015 in EUR bln. via robust coupling.**

| Year events | | 5 | 10 | 20 | 50 | 100 | 250 | 500 | 1000 |
|---|---|---|---|---|---|---|---|---|---|
| Cum. probabilities | | 0.8 | 0.9 | 0.95 | 0.98 | 0.99 | 0.996 | 0.998 | 0.999 |
| Fl. Clayton ($L$) | PNL=0.68 | 0.000 | 0.030 | 0.855 | 2.926 | 8.017 | 13.810 | 16.756 | 17.761 |
| Gumbel ($L$) | PNL=0.67 | 0.000 | 0.067 | 0.847 | 2.953 | 7.975 | 13.107 | 16.340 | 17.736 |
| Frank ($L$) | PNL=0.63 | 0.000 | 0.395 | 1.147 | 2.894 | 7.697 | 10.213 | 10.689 | 10.689 |

**Table 3** **Total losses in Austria for 2030 in EUR bln. via robust coupling.**

| Year events | | 5 | 10 | 20 | 50 | 100 | 250 | 500 | 1000 |
|---|---|---|---|---|---|---|---|---|---|
| Cum. probabilities | | 0.8 | 0.9 | 0.95 | 0.98 | 0.99 | 0.996 | 0.998 | 0.999 |
| Fl. Clayton ($L$) | PNL=0.68 | 0.000 | 0.028 | 0.850 | 3.004 | 8.123 | 13.875 | 16.853 | 17.877 |
| Gumbel ($L$) | PNL=0.67 | 0.000 | 0.063 | 0.844 | 3.039 | 8.074 | 13.173 | 16.435 | 17.850 |
| Frank ($L$) | PNL=0.63 | 0.000 | 0.406 | 1.138 | 2.967 | 7.781 | 10.300 | 10.785 | 10.785 |

**Table 4** **Total losses in Austria for 2050 in EUR bln. via robust coupling.**

| Year events | | 5 | 10 | 20 | 50 | 100 | 250 | 500 | 1000 |
|---|---|---|---|---|---|---|---|---|---|
| Cum. probabilities | | 0.8 | 0.9 | 0.95 | 0.98 | 0.99 | 0.996 | 0.998 | 0.999 |
| Year 2015 | PNL=0.62 | 17.400 | 21.000 | 24.100 | 27.400 | 29.600 | 32.400 | 34.400 | 34.400 |
| Year 2030 | PNL=0.62 | 17.700 | 21.300 | 24.400 | 27.800 | 30.100 | 32.900 | 35.000 | 35.000 |
| Year 2050 | PNL=0.62 | 17.600 | 21.200 | 24.300 | 27.900 | 30.200 | 32.100 | 35.200 | 35.200 |

**Table 5** **Total losses in Austria in EUR bln. in case of a strong dependency between extremes.**

| Year events | | 5 | 10 | 20 | 50 | 100 | 250 | 500 | 1000 |
|---|---|---|---|---|---|---|---|---|---|
| Cum. probabilities | | 0.8 | 0.9 | 0.95 | 0.98 | 0.99 | 0.996 | 0.998 | 0.999 |
| Year 2015 | PNL=0.60 | 0.000 | 0.433 | 1.140 | 2.830 | 7.190 | 9.710 | 10.500 | 10.500 |
| Year 2030 | PNL=0.59 | 0.000 | 0.419 | 1.150 | 2.820 | 7.130 | 10.100 | 10.700 | 10.700 |
| Year 2050 | PNL=0.59 | 0.000 | 0.406 | 1.120 | 2.890 | 7.220 | 9.790 | 10.800 | 10.800 |

**Table 6** **Total losses in Austria in EUR bln. in case of a weak dependency between extremes.**